\documentclass[11pt]{article}
\usepackage{graphicx} 
\usepackage{amsmath}
\usepackage{amsthm} 
\usepackage{amsfonts}
\usepackage[english]{babel}
\usepackage{xcolor}
\usepackage{gensymb}
\usepackage{hyperref}
 \usepackage{soul}

\newtheorem{defi}{Definition}
\newtheorem{prop}{Proposition}
\newtheorem{lem}{Lemma}

\def\cinfty#1{C^{\scriptscriptstyle\infty}(#1)}
\def\vectorfields#1{{\cal X}(#1)}

\def\fpd#1#2{{\displaystyle\frac{\partial #1}{\partial #2}}}

\def\sode{{\sc sode}}

\def\T{{\mathbf T}}
\def\D{{\mathcal D}}
\def\C{{\mathcal C}}
\def\GLD{{\Gamma_{(L,\mathcal{D})}}}

\parskip=\medskipamount
\def\I{I}
\def\J{J}
\def\m{M}
\def\Ra{R}
\def\ce{C_1}
\def\ct{C_2}
\def\cd{C_3}
\def\cv{C_4}

\begin{document}

\title{Geodesic extensions of mechanical systems with nonholonomic constraints}

\author{ Malika Belrhazi and Tom Mestdag \\[2mm]
	{\small Department of Mathematics,  University of Antwerp,}\\
	{\small Middelheimlaan 1, 2020 Antwerpen, Belgium}
}

\date{}

\maketitle

\begin{abstract} For a Lagrangian system with nonholonomic constraints, we {discuss} extensions of the equations of motion to  sets of second-order ordinary differential equations. In the case of a purely kinetic Lagrangian, we investigate the conditions under which the nonholonomic trajectories are geodesics of a Riemannian metric, while preserving the constrained Lagrangian. We interpret the algebraic and PDE conditions of this problem as infinitesimal versions of the relation between the nonholonomic exponential map and the Riemannian metric. We discuss  the special case of a Chaplygin system with symmetries and we end the paper with {some  worked-out examples}. 

\end{abstract}

\vspace{3mm}


\vspace{3mm}


\maketitle

\section{Introduction}

 In many applications involving mechanical systems, the motion is  often limited by velocity-dependent constraints. If these are given by equations, and if they are non-integrable, they are called nonholonomic constraints. Some examples of such systems are those that involve wheels that roll without slipping, or a skate that is prevented to move in a direction perpendicular to the blade. 

Let $L=T-V$ be the Lagrangian of such a mechanical system (given by the difference in kinetic and potential energy) on a configuration manifold $Q$ and let  the  nonholonomic constraints (assumed to be linear, throughout the paper) be represented by a distribution $\D$ on $Q$.  The equations of motion for such systems $(L,\D)$ follow from the so-called Lagrange-d'Alembert principle (see e.g. \cite{Bloch,Cortes}). It is well-known that this is not a purely variational principle (for the Lagrangian functional associated to $L$), since it only makes use of  variations that are required to follow the constraints. The resulting equations, the so-called Lagrange-d'Alembert equations of $(L,\D)$, are therefore not simply the Euler-Lagrange equations of $L$. Instead, the nonholonomic equations of motion contain Lagrangian multipliers, and when eliminated, they form a set of first- and second-order ordinary differential equations (see their expressions in Section~\ref{sec2}).

We will call  the solutions of the Lagrange-d'Alembert equations the nonholonomic trajectories of $(L,\D)$. In many papers, it has been recognized that it is often useful if the nonholonomic trajectories can be interpreted, in one way or another, as part of  the solutions of the Euler-Lagrange equations of some unrelated (and possibly non-mechanical) Lagrangian $\hat L$. This question is often referred to as the `Hamiltonization problem of nonholonomic mechanics'. It was probably first raised in the context of quantization of classical nonholonomic systems (see \cite{abud,BR,quant}), but there are also other important areas of  possible applications, such as e.g.\ that of geometric integrators of nonholonomic systems (see e.g. \cite{FernandezOlver,TomBrussel}). 
The subtlety of the Hamiltonization problem  lies in the fact that, in contrast to the nonholonomic equations of $(L,\D)$, the Euler-Lagrange equations of $\hat L$ are purely second-order ordinary differential equations. 
There have been many different approaches to this question. For example, in \cite{BFM,marta} the main idea is to enlarge the mixed first- and second order system of $(L,\D)$  to a (chosen) full system of \sode s and to find a Lagrangian $\hat L$ for that system  using techniques from the so-called  inverse problem of Lagrangian mechanics. In other papers \cite{BalseiroFernandez,BalseiroNaranjo,Bol,Ehlers} one rather uses the symmetries of the nonholonomic system to (roughly speaking) first  reduce it to  a \sode\ on a quotient manifold. Then, one tries to find a symplectic or Poisson structure on that quotient for which it is a Hamiltonian vector field. 

The results of the current paper belong to the first approach. In short, we will investigate the conditions under which it is possible to have a Hamiltonization by means of the geodesic equations of a Riemannian metric. 
If that can be achieved,  we will speak of  a `geodesic extension' of the nonholonomic dynamics (even though this question seems to be dubbed the kinetic Lagrangianization  problem in the conclusions of the paper \cite{Simoes}).

In Section~\ref{sec2} we  recall  the most elementary concepts in this context, in an attempt to keep the  paper self-contained. Here it is important to notice that we will use a specific technique throughout the paper. Indeed, to suppress the use of long expressions, and to gain ready geometric insight, it will be advantageous to express the  intrinsic and coordinate-independent tools that we will use in terms of a so-called anholonomic frame of vector fields. In particular, we will often choose the frame in such a way that the vanishing of a part of its corresponding quasi-velocities represents the constraints and its distribution $\D$. In the last sections, it will be of further advantage to specify the frame in such a way that it also follows the orthogonal decomposition of the tangent spaces by a given Riemannian metric.

In Section~\ref{sec3} we will first define the notion of a \sode\ extension  of a nonholonomic system as a set of (only) second-order differential equations whose solutions set contains the nonholonomic trajectories. The name comes from \cite{Simoes}, but the idea can already be found in e.g.\ \cite{BC}  as `Newton laws with the restriction property' and, implicitly, in \cite{BFM} where the \sode\ extensions are called `associated systems'.  {There are many more places in the literature where they appear in a somewhat disguised form, such as in e.g.\ \cite{Fasso}}. In the first part of the paper, we will show that such \sode\ extensions always exist. We will give two natural constructions of a \sode\ extension in Propositions~\ref{prop:1}  and \ref{prop2}. We will motivate our constructions by showing in Proposition~\ref{prop:so} that they represent, in fact, a generalization of  the  well-known nonholonomic connection (see e.g.\cite{Lewis,Cortes}) to the level of arbitrary Lagrangians. We will say that a nonholonomic system is purely kinetic, if its potential $V$ is  constant and if its kinetic energy can be derived from a Riemannian metric $g$. Specifically for a  purely kinetic nonholonomic systems, there is also a third construction of a \sode\ extension at the end of  Section~\ref{sec6}.

From Section~\ref{sec4}  onwards, we restrict the attention to only  purely kinetic nonholonomic systems.  We will investigate the conditions under which the Riemannian metric $g$ of the kinetic energy can be modified to a metric $\hat g$, such that the nonholonomic trajectories become part of its set of  geodesics (and therefore also of the Euler-Lagrange equations that can be associated to the Lagrangian $\hat L$ of the kinetic energy $\hat g$). The specific kind of geodesic extension we envisage here will be such that the  constrained Lagrangian $L_C = L|_\D$ remains preserved. This means that we will assume that the restrictions of the Riemannian metrics $g$ and $\hat g$ to the constraint distribution $\D$ are the same. We will say, in that case, that $\hat g$ is a $\D$-preserving modification of $g$.
  We are able to show in Lemma~\ref{prop:result} that this has the advantage that the question narrows  to the search for a suitable 1-form  (with components $\theta_i$ along the distribution $\D$) that satisfies a set of algebraic and PDE conditions. If these have a solution, we show in Proposition~\ref{prop:result2}  that we can always find  a Riemannian metric $\hat g$ that satisfies our purposes.  At the end of the section we indicate how one can generate more necessary algebraic conditions.  In Section~\ref{sec5} we relate (in Propositions~\ref{propback1} and \ref{propback2}) the conditions of Lemma~\ref{prop:result} to the kind of  \sode\ extensions  we had encountered in Section~\ref{sec3}.

 The exponential map of a Riemannian metric  plays an important role in differential geometry and geometric mechanics, for example in results concerning Jacobi fields, conjugate points, variations of the energy and focal points. In this paper we will also address the search for a geodesic extension from the viewpoint of the so-called nonholonomic exponential map. This concept has been recently introduced in the paper \cite{firstnonholexp}, and further developed in \cite{Simoes}. There, among other results, it was used to show that the nonholonomic trajectories of a purely kinetic system, originating from a fixed point $q$, can be interpreted as the radial geodesics of a Riemannian metric on a submanifold. This Riemannian metric should not be confused with the (full) Riemannian metric $\hat g$ we are looking for, since it (and the submanifold on which it is defined) is different for each different point $q$. Notwithstanding, like \cite{Simoes}, we will be inspired here by the same underlying result from Riemannian geometry: the so-called Gauss Lemma (see e.g. \cite{do}). This Lemma states that the exponential map of a Riemannian metric at $q$  is  a radial isometry. In Proposition~\ref{prop:gauss} we show that the algebraic condition of Lemma~\ref{prop:result} follows, as soon as the exponential map and the sought-for modified metric $\hat g$ satisfy a certain Gauss-type condition.  Next, in Proposition~\ref{prop:PT} we show that also the PDE condition can be interpreted as an infinitesimal condition. Along each nonholonomic trajectory, it can be seen to be equivalent with a conservation law that involves the parallel transport along a linear connection whose quadratic spray is, in fact, again a certain \sode\ extension.  For these reasons, we will call the algebraic and PDE conditions   infinitesimal Gauss conditions.

In Section~\ref{sec7} we investigate the presence of symmetry. Chaplygin nonholonomic systems are such that the Lagrangian and the constraints are both invariant under the action of a symmetry Lie group $G$. Moreover, it is assumed that $Q\to Q/G$ is a principal fibre bundle and that the constraint distribution $\D$ is the horizontal distribution of a principal connection. We show in Proposition~\ref{propchap} that the conditions of Lemma~\ref{prop:result} can be simplified, and that the geodesic extension problem can be stated as the search for linear first integrals that satisfy some further algebraic conditions.

We illustrate the results of the paper with the example of a disk that rolls vertically without slipping, and we compare the metrics we find with those that had already appeared in the literature for this example \cite{Simoes,recent,FB,BFM}.  With the example of the two-wheeled carriage we show how we can benefit from the extra necessary algebraic conditions to find an answer to the geodesic extension question. We end the paper with a discussion on some of the limitations of our approach, and we give a few directions for future research.

\section{Tangent bundle geometry,  Lagrangian systems and nonholonomic systems} \label{sec2}

We will assume throughout  that $Q$ is an $n$-dimensional real differentiable  manifold, with local coordinates $(q^\alpha)$. We will assume that $\{X_\alpha\}$ represents an  anholonomic frame of vector fields  on $Q$. The adjective anholonomic means here that the Lie brackets of these vector fields do not necessarily vanish, and satisfy the relation $[X_\alpha,X_\beta] = R^\gamma_{\alpha\beta}X_\gamma$. We will call the components $v^\alpha$ in the decomposition $v_q = v^\alpha X_\alpha (q)$ the quasi-velocities of a tangent vector $v_q\in T_qQ$. We may then use $(q^\alpha,v^\alpha)$ as (non-standard) coordinates on the tangent manifold $TQ$, where most of the objects of interest live.

Any vector field $X$ on $Q$ can be lifted to vector fields on $TQ$, by means of the complete and vertical lift, respectively: If  $X=X^\alpha\partial/\partial{q^\alpha}$, then in natural coordinates $(q^\alpha, {\dot q}^\alpha)$,
$$
 X^C=X^\alpha \frac{\partial}{\partial q^\alpha}+\dot{q}^\beta \frac{\partial X^\alpha}{\partial q^\beta}\frac{\partial}{\partial \dot{q}^\alpha} \quad\text{and}\quad X^V=X^\alpha \frac{\partial}{\partial \dot{q}^\alpha}.
$$ 
{In the following, we will assume that $\{X_\alpha\}$ is a frame for $\mathcal{X}(Q)$ (the set of vector fields on $Q$), in the sense that each vector field $X\in\mathcal{X}(Q)$ can be decomposed as $X=\xi^\alpha X_\alpha$, for some smooth functions $\xi^\alpha$ on $Q$.} One may easily verify that, if $\{X_\alpha\}$ is a {frame} for $\mathcal{X}(Q)$, then $\{X_\alpha^C,X_\alpha^V\}$ forms a   {frame} for $\mathcal{X}(TQ)$. An important property we will frequently use is that, on the quasi-velocities $v^\alpha$ of the frame,
$$
X_\alpha^V(v^\beta)=\delta_\alpha^\beta,\quad\quad\quad X_\alpha^C(v^\beta)=-R_{\alpha\gamma}^\beta v^\gamma,
$$
see e.g.\ \cite{anhol} (and the references therein)  for more details.

We recall some aspects of the calculus of forms along the tangent bundle projection $\tau:TQ\rightarrow Q$ (see e.g. \cite{derII,Willysurvey} for a survey). 
A section $X$ of the pullback bundle $\tau^* TQ=\{(u,v)\in TQ\times TQ\,|\, \tau(u)=\tau(v)\} \to TQ$ is said to be a vector field along $\tau$. We denote the set of all these vector fields by $\mathcal{X}(\tau)$. We may interpret any vector field $Y$ on $Q$ as a ``basic'' vector field $Y\circ\tau$ along $\tau$, but we will no longer make a notational distinction. For that reason, we may express $X\in \vectorfields{\tau}$ in the given frame $\{X_\alpha\}$ as  
$X=\xi^\alpha X_\alpha$, with   now  $\xi^\alpha \in \cinfty{TQ}$.
It is easy to see that vertical vector fields on $TQ$ are in 1-1 correspondence with vector fields along $\tau$. Indeed, each $X=\xi^\alpha X_\alpha\in\mathcal{X}(\tau)$ can be vertically lifted to $X^V\in\mathcal{X}(TQ)$, given by 
$$
X^V=\xi^\alpha X_\alpha^V.
$$ 
In particular, there exist a canonical vector field $\T = v^\alpha X_\alpha\in\mathcal{X}(\tau)$. It  corresponds to the section $v_q\mapsto (v_q,v_q)$ of the pullback bundle, and its   vertical lift  gives the Liouville vector field on $TQ$.

  A vector field $\Gamma$ on $TQ$ is a second-order ordinary differential equations field (\sode, in short) if all its integral curves $\gamma: I\rightarrow TQ$ are lifted curves, that is of the type $\gamma=\dot{c}$. We   call  $c:I\rightarrow Q$ a base integral curve of $\Gamma$.  
When expressed  in terms of the frame $\{X_\alpha^C,X_\alpha^V\}$, a \sode\ $\Gamma$ is of the form 
$$
\Gamma=v^\alpha X_\alpha^C+F^\alpha X_\alpha^V,
$$
for some functions $F^\alpha\in\mathcal{C}^\infty(TQ)$. The quasi-velocity components $(v^\alpha(t))$ of the integral curves $\gamma(t)$ of $\Gamma$ satisfy
$\dot{v}^\alpha=F^\alpha $.

A \sode\ $\Gamma$ can be used to introduce a horizontal lift $X^H\in\mathcal{X}(TQ)$ of  $X\in\mathcal{X}(\tau)$ (see e.g. \cite{derII,Willysurvey}). Its coordinate expression in the anholonomic frame is
$$
X^H=\xi^\alpha \left(X_\alpha^C-\Gamma^\beta_\alpha X_\beta^V \right),\quad\quad \Gamma^\beta_\alpha=-\frac{1}{2} \Big( X_\alpha^V(F^\beta) + v^\gamma R^\beta_{\alpha\gamma}\Big) .
$$
With this, any vector field $Z$ on $TQ$ can be decomposed into the sum of a horizontal and a vertical lift,  $Z=X_1^H+X_2^V$ for $X_1,X_2\in\mathcal{X}(\tau)$. For any vector field $X$ along $\tau$, we may, simultaneously, define  two operators $\nabla$ and $\Phi$ by the relations  
$$
[\Gamma, X^V]=-X^H+(\nabla X)^V\quad\text{and}\quad [\Gamma,X^H]=(\nabla X)^H+(\Phi(X))^V.
$$
 The operator $\Phi$ is a  $(1,1)$-tensor field along $\tau$, the so-called Jacobi endomorphism. It is of no concern for the current paper. Second, from the above defining relation of the dynamical covariant derivative $\nabla$ of $\Gamma$, one can derive that it acts  as a derivation on $\mathcal{X}(\tau)$. This means that, for $f\in\cinfty{TQ}$ and $X\in\vectorfields{\tau}$:
$$
\nabla(fX)=f\nabla X+\Gamma(f)X\quad\quad\text{and}\quad\quad \nabla X_\alpha=\Gamma^\beta_\alpha X_\beta.
$$

The equations of motion for an (unconstrained) Lagrangian systems (with  regular Lagrangian $L:TQ\rightarrow \mathbb{R}$) are given by the Euler-Lagrange  equations:   
$$
\frac{d}{dt}\left(\frac{\partial L}{\partial \dot{q}^\alpha}\right)-\frac{\partial L}{\partial q^\alpha}=0.
$$ 
The solutions of these equations are base integral curves of the Lagrangian vector field $\Gamma_L\in\mathcal{X}(TQ)$, which can be defined by as the unique \sode\ that satisfies
$$
\Gamma_L(X^V(L))-X^C(L)=0,\quad\forall X \in\mathcal{X}(Q).
$$

The main focus of the paper is on nonholonomic mechanical systems.
\begin{defi}
 A nonholonomic system on a manifold $Q$ consists of a pair $(L,\mathcal{D})$, where $L:TQ\rightarrow \mathbb{R}$ is the Lagrangian of the system and $\mathcal{D}$ is the $m$-dimensional distribution on $TQ$ that describes the nonholonomic constraints.  
\end{defi}
In what follows, we will often use $\C$, when we interpret the distribution $\D$ specifically as a submanifold of $TQ$. We will also use $k$ for  $n-m$.

 The equations of motion follow from considering an extended principle of Hamilton (also called Lagrange-d'Alembert principle, see \cite{Bloch,Cortes}). The nonholonomic trajectory  $c_{v_q}(t) = (q^\alpha(t))$  that starts at $v_q\in \D_q$ is a solution of the equations 
\begin{equation}
\begin{cases}
\displaystyle\mu_\alpha^i\dot{q}^\alpha=0,
\\
\displaystyle\frac{d}{dt}\Big(\frac{\partial L}{\partial \dot{q}^\alpha}\Big)-\frac{\partial L}{\partial q^\alpha}=\sum_{i=1}^k{\Lambda_i} \mu^i_{\alpha}.
\end{cases}\nonumber 
\end{equation}
The first set of equations represent the nonholonomic constraints. They ensure that the velocity ${\dot c}_{v_q}(t)$ remains in the distribution $\D$ during the motion. In the second set of equations, {the functions $\Lambda_i(t)=\Lambda_i(q(t),\dot{q}(t))$ }are called  the Lagrange multipliers. 

As before for Lagrangian systems, we will represent the nonholonomic trajectories  as the base integral curves of a vector field $\Gamma_{(L,\mathcal{D})}$. Suppose that the vector fields $\{X_a\}$ span $\mathcal{D}$, and that we complete this set to a frame $\{X_\alpha\} = \{X_a,X_i\}$ for all vector fields on $Q$. This has the advantage that the equations for the nonholonomic constraints become very simple: if $(v^\alpha)=(v^a,v^i)$ are the quasi-velocities with respect to this frame, then $v_q\in D_q$ if and only its quasi-velocities satisfy  $v^i=0$. In this paper, the indices $\alpha$ run from 1 to $n$, $a$ from 1 to $m$ and $i$ from 1 to $k$.

We will say that the Lagrangian is regular with respect to $\D$ if the matrix
\[
\left(X^V_a(X^V_b(L))\right)
\]
is non-singular everywhere on $TQ$. This is actually a small modification from the terminology in \cite{anhol,Cartan}, where the matrix is only assumed to be non-singular on $\C$ (and see also the regularity requirements in \cite{deL,Cortes}). If satisfied, there will exist a unique vector field $\GLD$ of the type
$$
\Gamma_{(L,\mathcal{D})}=v^{a} X_{a}^C+f^{a} X_{a}^{V}
$$
that satisfies, on $\C$,  
$$
\Gamma_{(L,\mathcal{D})}(X^V(L))-X^C(L)=0 ,\quad\forall X\in{\rm Sec}(\mathcal{D}).
$$
We will call this vector field the nonholonomic vector field. Given the current frame, it can be completely determined from the relation
$$
\Gamma_{(L,\mathcal{D})}(X^V_a(L))-X^C_a(L)=0, \qquad \mbox{on $\C$}.
$$ 

This equation represents a version of the Lagrange-d'Alembert equations where the Lagrangian multipliers have already been eliminated. Once $\GLD$ is determined, the functions $\lambda_i:=\Gamma_{(L,\mathcal{D})}(X^V_i(L))-X^C_i(L)$ on $\C$ play the role of the Lagrangian multipliers, in the current frame.

Let $L$ be a regular Lagrangian. The Hessian $g$ of $L$, with respect to fibre coordinates, can be thought of as a $(0,2)$-tensor field along the tangent bundle projection $\tau:TQ\rightarrow Q$ (see \cite{anhol,Willysurvey} for details). To define it, it is enough to state its action along `basic' vector fields in $\vectorfields{\tau}$ (i.e.\ vector fields on $Q$): 
\[
g(X,Y) = X^V(Y^V(L)), \qquad X,Y\in\vectorfields{Q}.
\]
and extend it to $\vectorfields{\tau}$ by $\cinfty{TQ}$-linearity.
Since there is a 1-1 correspondence between $\tau$-vertical vector fields on $TQ$ and vector fields along $\tau$, we may also think of it as acting on vertical vector fields. 

For what follows, it is important to realize that this Hessian is familiar in the case of a  purely kinetic Lagrangian, i.e.\ when $L(v)=\frac{1}{2}G(v,v)$ for a Riemannian metric $G$ on $Q$. In this case, we may identify the above defined Hessian $g$ of $L$ (interpreted as a (0,2) tensor field along $\tau$) with the `basic' Riemannian metric $G$ (interpreted as a (0,2) tensor field on $Q$). We will not make a notational distinction between $G$ and $g$ in that case, and simply write  $L(v)=\frac{1}{2}g(v,v)$.

Let's assume now that $L$ is again  an arbitrary regular Lagrangian that is also regular with respect to $\D$. Let $g$ be its Hessian. The pullback bundle $\tau^*\D \to TQ$ is a vector subbundle of $\tau^*TQ$ whose set of sections we will denote by $\vectorfields{\tau,\D}$. Since $L$ is regular with respect to $\D$, we can consider the set  ${\mathcal X}^g(\tau,\D)$  of vector fields along $\tau$ that are orthogonal to $\vectorfields{\tau,\D}$, with respect to $g$. That is, $X\in {\mathcal X}^g(\tau,\D)$ when
\[
g(X,Y)=0,\qquad \forall Y \in \vectorfields{\tau,\D}.
\]   
With this, we can decompose $\vectorfields{\tau}$ as
\[
\vectorfields{\tau} = \vectorfields{\tau,\D} \oplus {\mathcal X}^g(\tau,\D)
\]
and define two complementary projection operators: 
\[P:\vectorfields{\tau} \to \vectorfields{\tau,\D}\quad\mbox{ and }\quad Q: \vectorfields{\tau} \to {\mathcal X}^g(\tau,\D).
\]

Starting from a frame $\{X_a,X_i\}$ for $\vectorfields{Q}$, we  can easily construct a  frame for $\vectorfields{\tau}$ that is adapted to the decomposition. Indeed, if $g_{ab}=g(X_a,X_b)$, $g_{ai} = g(X_a,X_i)$ and $g_{ij} =g(X_i,X_j)$, then $g_{ab}$ has an inverse (because of the assumed regularity with respect to $\D$). Then, the set  $\{X_a, {\tilde X}_i\}$, with 
\[
\tilde{X}_i=X_i- {g}^{a b} {g}_{b i} X_a = X_i+ K^a_i X_a
\]
is such a frame. The vector fields $\tilde{X}_i\in\vectorfields{\tau}$ are, however, not basic vector fields.

Next, we introduce the set ${\mathcal X}^T(TQ)$ of vector fields  on $TQ$ that are tangent to the constraint manifold $\C$, i.e.\ 
$$
{\mathcal X}^T(TQ)=\{Z\in\mathcal{X}(TQ)\,|\, Z(v^i)=0\}.
$$

In view of the relations  $X_\alpha^C(v^i)=-R^i_{\alpha\gamma}v^\gamma$ and $X_\alpha^V(v^i)=\delta^i_\alpha$, one easily verifies that the vector fields on $TQ$, given by
\[
    \overline{X_i^C}= X_i^C+R^j_{i\gamma}v^\gamma \tilde{X}_j^V,\quad
    \overline{X_a^C}= X_a^C+R^j_{a\gamma}v^\gamma \tilde{X}_j^V \quad \mbox{and}\quad  X_a^V
\]
form a basis for ${\mathcal X}^T(TQ)$. 
We may augment this basis with the vector fields $\{{\tilde X_i}^V\}$ that span $({\mathcal X}^g(\tau,\D))^V$ to obtain a full basis of $\vectorfields{TQ}$. From this, we can conclude that we can decompose $\mathcal{X}(TQ)$ as
   $$
	{\mathcal X}(TQ)={\mathcal X}^T(TQ)\oplus ({\mathcal X}^g(\tau,\D))^V
	$$
and 
that we can define two complementary projectors $p$ and $q$,
$$p:{\mathcal X}(TQ)\rightarrow {\mathcal X}^T(TQ)\quad\quad \text{and}\quad\quad q: \mathcal{X}(TQ)\rightarrow ({\mathcal X}^g(\tau,\D))^V.$$
If $Z=Z^a X_a^C+Z^i X_i^C+\zeta^a X_a^V+\zeta^i X_i^V  \in \mathcal{X}(TQ)$, then 
\[
p(Z) = Z^a\overline{X_a^C} +Z^i\overline{X_i^C}+ (\zeta^a-\zeta^i K_i^a) X_a^V\quad \mbox{and} \quad q(Z)= (-Z^a R^j_{a\gamma}v^\gamma-Z^i R^j_{i\gamma}v^\gamma+\zeta^j)\tilde{X}_j^V.
\]

\section{Extensions of the nonholonomic vector field to a \sode} \label{sec3}

In what follows, we will  assume throughout that  the Lagrangian $L$ is both regular and regular with respect to $\D$.  
Let $c_{v_q}(t)$ be the nonholonomic trajectory that starts at $v_q\in \D_q$.  We formalize the notion of a \sode\ extension (implicit in  \cite{firstnonholexp,Simoes})  in the following  definition.

\begin{defi} A \sode\ $\Gamma$   is a \sode\ extension of the nonholonomic vector field $\Gamma_{(L,\mathcal{D})}$, if each nonholonomic trajectory $c_{v_q(t)}$ is a base integral curve of $\Gamma$.\label{defi:dfn}
\end{defi}

An equivalent way to define a \sode\ extension is that it should satisfy
\[
\Gamma|_\mathcal{D}=\Gamma_{(L,\mathcal{D})}.
\]
It's   useful to see this property translated into terms of a frame  $\{X_\alpha\}=\{X_a,X_i\}$ on $Q$ whose first elements $\{X_a\}$ span $\D$, and whose corresponding quasi-velocities are $(v^\alpha)=(v^a,v^i)$. The \sode\ $\Gamma$ and the vector field $\Gamma_{(L,\mathcal{D})}$  can then be represented, respectively, by
\[\Gamma=v^\alpha X_\alpha^C+F^\alpha X_\alpha^V \quad\mbox{and}\quad \Gamma_{(L,\mathcal{D})}=v^a X_a^C+ f^a X_a^V.
\]
Let $v^a(t)$ be such that $\dot{v}^a=f^a$. From the definition, $\Gamma$ is a \sode\ extension   if  the curve $(v^a(t),v^i(t)=0)$ is an integral curve of $\Gamma$, i.e.\ a solution of  $\dot{v}^a=F^a$ and $\dot{v}^i=F^i$. This will be the case, if and only if, 
\[
    F^a|_\mathcal{C}=f^a,\qquad \mbox{and} \qquad
    F^i|_\mathcal{C}=0.
 \]

For the rest of the paper it is important to realize that the Lagrangian vector field $\Gamma_L$ of $L$ is, {in general}, {\em not} a \sode\ extension of $\Gamma_{(L,\mathcal{D})}$. The \sode\ $\Gamma_L =v^\alpha X_\alpha^C+F^\alpha X_\alpha^V$ is uniquely determined by the relations
$${\Gamma_{L}\left(X_a^V(L)\right)-X_a^C(L)=0,\qquad \Gamma_{L}\left(X_i^V(L)\right)-X_i^C(L)=0.}
$$
From the first set of equations we get, in particular  on $\C$,
$$v^bX_b^C (X_a^V(L))+F^b|_\mathcal{C} X_b^V (X_a^V(L))+ F^i|_\mathcal{C}X_i^V (X_a^V(L))-X_a^C(L)=0. 
$$
On the other hand, on $\C$, the nonholonomic vector field $\Gamma_{(L, \mathcal{D})}= v^aX_a^C+f^aX_a^V$ satisfies:
$$
\Gamma_{(L, \mathcal{D})}\left(X_a^V(L)\right)-X_a^C(L)=0,
$$
which becomes when written in full
$$
v^bX_b^C (X_a^V(L))+ f^b X_b^V (X_a^V(L))-X_a^C(L)=0.
$$
In the assumed regularity with respect to $\D$, this last equation determines the coefficients $f^b$ of the nonholonomic vector field. This means that the relation between the forces $F^a$, $F^i$ and $f^a$ is 
\[
f^a= F^a|_\C-K^a_i F^i|_\C \quad \text {on } \mathcal{C}.
\]
This is not the required property for $\Gamma_L$ to be a \sode\ extension of $\GLD$, but it is a useful property for what follows. {Although $\Gamma_L$ is not a \sode\ extension in general, it may still happen that it is one, in some specific examples.}

We now show that a \sode\ extension of $\GLD$ always exists.
The construction we present below is, in essence, also contained in  \cite{Cartan} (Proposition~2) and \cite{deL} (formulated as Proposition~2.1 in \cite{firstnonholexp}).  
 
\begin{prop} \label{prop:1} The vector field $\Gamma^1 =p(\Gamma_L)$ is a \sode\ extension of the nonholonomic vector field $\Gamma_{(L,\mathcal{D})}$. It can, equivalently, be defined as the unique \sode\ $\Gamma^1$ which is tangent to the constraints, and whose difference from $\Gamma_L$ is $g$-orthogonal. 
\end{prop}

\begin{proof}
If we rewrite $\Gamma_L$ in terms of the basis $\{\tilde{X}_i^V, \overline{X_i^C},\overline{X_a^C}, X_a^V\}$ we  get 
\begin{eqnarray*}
    \Gamma^1&=& p(\Gamma_L)= v^a\overline{X_a^C} +v^i\overline{X_i^C}+ (F^a-F^i K_i^a) X_a^V\\
    &=& v^a(X_a^C+R^j_{a\gamma}v^\gamma \tilde{X}_j^V) +v^i(X_i^C+R^j_{i\gamma}v^\gamma \tilde{X}_j^V)+ (F^a-F^i K_i^a) X_a^V\\
    &=& v^a X_a^C+v^i X_i^C+v^\mu R^j_{\mu\gamma}v^\gamma \tilde{X}_j^V+(F^a-F^i K^a_i)X_a^V\\
    &=& v^\alpha X_\alpha^C+(F^a-F^i K_i^a)X_a^V+0X_i^V.
\end{eqnarray*}
From the first term we may conclude that $\Gamma^1$ is indeed a \sode. When restricted to $\C$, we know that $F^a-F^i K_i^a$ coincides with $f^a$. Therefore, we can conclude that $p(\Gamma_L)$ is a \sode\ extension.

To prove the second statement, we set $\Gamma^1=\Gamma_L+U^V=v^\alpha X_\alpha^C+F^\alpha X_\alpha^V+U^\alpha X_\alpha^V$. From the first condition, $\Gamma^1(v^i)=0$, we obtain $U^i=-F^i$. The second condition can be  expressed as  $U\in({\mathcal X}^g(\tau,\D))^V$, and thus $ U^b=K^b_iU^i = - K^b_iF^i$. With this we see that, again,
$
\Gamma^1= v^\alpha X_\alpha^C+(F^a-F^i K_i^a)X_a^V.
$
\end{proof}

In the next paragraph, we give a second construction of a \sode\ extension and show how it appears in the literature (in a special case).

\begin{prop}\label{prop2}
Let $\nabla$ be the dynamical covariant derivative of $\Gamma_L$ and $Q$ the projection on ${\mathcal X}^g(\tau,\D)$, defined by $g$. The vector field 
$$
\Gamma^2=\Gamma_L-[\nabla(Q(\T))]^V
$$
is a \sode\ extension of the nonholonomic vector field $\Gamma_{(L,\mathcal{D})}$.
\end{prop}

\begin{proof}
For the  canonical vector field $\T$ along $\tau$ we can write
$
\T= v^a X_a+v^i X_i = \tilde{v}^a X_a+\tilde{v}^i \tilde{X}_i
$
with $\tilde{v}^a=v^a+v^i  K^a_i$ and $\tilde{v}^i=v^i$. Remark that if $v^i=0$ (on $\mathcal{C}$) then $\tilde{v}^i=0$ and $\tilde{v}^a=v^a$. From this, $Q(\T) = v^i {\tilde X}_i$ and
\[
\nabla(Q(\T))=\nabla(v^i \tilde{X}_i) = \Gamma_L(v^i) \tilde{X}_i+v^i \nabla \tilde{X}_i= F^i \tilde{X}_i+v^i\nabla \tilde{X}_i.
\]
The last term doesn't really concern us,  since it vanishes on $\C$. Then,
\[
\Gamma^2 =v^\alpha X_\alpha^C+F^a X_a^V+F^i X_i^V-F^i \tilde{X}_i^V -v^i(\nabla \tilde{X}_i)^V =v^\alpha X_\alpha^C+(F^a -F^i K^a_i) X_a^V-v^i(\nabla \tilde{X}_i)^V,
\]
from which  it is clear that $\Gamma^2|_{\C}=\Gamma_{L,\D}$.
\end{proof}

The above construction is actually familiar, in the case of a  purely kinetic Lagrangian.

\begin{defi}
A nonholonomic system $(L,\D)$ is purely kinetic if its Lagrangian is of the type $L(v)=\frac{1}{2}g(v,v)$, for a Riemannian metric $g$ on $Q$.
\end{defi}
We have already remarked before that, in this case, the Hessian of $L$ (interpreted as a (0,2) tensor field along $\tau$) is directly related to the Riemannian metric $g$. Since $g$ is a Riemannian metric, and therefore positive-definite, $L$ is both regular and regular with respect to $\D$.  We first derive an expression for $\Gamma_L$ and $\GLD$. Next, we introduce the nonholonomic connection $\nabla^{nh}$, and we show that the \sode\ extension $\Gamma^2$ is in fact the quadratic spray of  $\nabla^{nh}$. 

We may use the Riemannian metric $g$ to define an orthogonal complement $\D^g$ of $\D$ in $TQ$. With a slight abuse of notation, we will write $P:\vectorfields{Q}\to {\rm Sec}(\D)$ and $Q:\vectorfields{Q}\to {\rm Sec}(\D^g)$ for the two corresponding projections, since these operators can clearly be extended to their counterparts we had encountered in the previous section (with the same notation) on $\vectorfields{\tau}$.

Many papers make use of the so-called  nonholonomic connection $\nabla^{nh}$ (see e.g. \cite{Lewis,Cortes}). It can be  defined by the expression 
$$
\nabla_X^{nh} Y=P\left(\nabla_X^g Y\right)+\nabla_X^g(Q(Y)).   
$$
Herein, $\nabla^g$ is the Levi-Civita connection of $g$.

We will assume (until the end of the paper)  that the frame $\left\{X_{\alpha}\right\} = \left\{X_{a},X_i\right\}$ is chosen is such a way that $\left\{X_{a}\right\}$ represents a frame for $\mathcal{D}$ and that $\left\{X_{i}\right\}$ is a frame for $\mathcal{D}^g$. This means that, from now on,  $g_{a i}=g\left(X_{a}, X_{i}\right)=0$.
As before, we denote the quasi-velocities w.r.t.\ this frame by $(v^\alpha)=(v^a,v^i)$.

It is well-known that the  Levi-Civita connection of $g$ satisfies the Koszul-formula
\[{2 g(\nabla_{X}^{g} Y, Z)=X(g (Y, Z))+Y(g(X, Z))-Z (g (X, Y))+g([X, Y], Z)  +g([Z, X], Y)+g([Z, Y], X) .}
\]
If we write $\nabla^g_{X_a}X_b = \Gamma_{ab}^cX_c + \Gamma_{ab}^i X_i$, $[X_a,X_b] =R_{ab}^c X_c+ R_{ab}^i X_i$, and so forth, we get by taking $X=X_a, Y=X_b, Z=X_c$ (among other choices) 
$$
2 g_{c d} \Gamma_{a b}^{d}=X_{a}(g_{b c})+X_{b}(g_{a c})-X_{c}(g_{a b})
+g_{d c} R_{a b}^{d} + g_{d b} R_{c a}^{d} + g_{d a} R_{c b}^{d}. 
$$
After multiplying both sides by $v^av^b$ and canceling a factor 2, we obtain
$$g_{c d} \Gamma_{a b}^{d} v^{a} v^{b}=X_{a}\left(g_{b c}\right) v^{a} v^{b}-\frac{1}{2} X_{c}\left(g_{a b}\right) v^{a} v^{b} 
-g_{d b} R_{a c}^{d} v^{a} v^{b}.$$

\begin{lem} \label{lem:VF} For a  purely kinetic Lagrangian
$
L=\frac{1}{2} g_{\alpha\beta}v^\alpha v^\beta = \frac{1}{2}\left(g_{ij} v^{i} v^{j}+g_{a b} v^{a} v^{b}\right) 
$
the  Lagrangian vector field is the geodesic spray $\Gamma_g$ of the Levi-Civita connection. It is given by the expression 
\[
\Gamma_L= \Gamma_g= v^\alpha X^C_\alpha - \Gamma^\alpha_{\beta\gamma}v^\beta v^\gamma X^V_\alpha.
\]
The nonholonomic vector field is given by
 $$\Gamma_{(L, \mathcal{D})}=v^{d} X_{d}^{C}-\Gamma_{a b}^{d} v^{a} v^{b}  X_{d}^{V}.
$$ 
The spray of the nonholonomic connection $\nabla^{nh}$ is of the form 
$$\Gamma_{\nabla^{nh}}=v^\gamma X_\gamma^C-\left(\Gamma_{a b}^{c} v^a v^b+2 \Gamma_{a j}^c v^a v^j+\Gamma_{i b}^c v^i v^b+2 \Gamma_{ij}^c v^iv^j\right) X_c^V 
-\left(\Gamma_{a j}^k v^a v^j+\Gamma_{ij}^k v^i v^j\right) X_k^V.
$$ 
\end{lem}

\begin{proof} We only prove the second expression, since the expression for the first can be found by a similar reasoning. 

 Recall that on $\mathcal{C}$  all $v^i=0$ and that 
$X_d^C(v^b)=-R_{de}^b v^e$ and $X_d^V(v^b)=\delta^b_d$. 
If we set $\Gamma_{(L_g, \mathcal{D})}=v^{d} X_{d}^{C}+f^d  X_{d}^{V}$, the coefficients $f^d$ satisfy on $\mathcal{C}$
\begin{eqnarray*}
0&=& \Gamma_{(L_g,\mathcal{D})}\left(X_{a}^{V}(L)\right)-X_{a}^{C}(L)\\
&=& v^{d} X^C_{d} (g_{a b} v^{b})+f^d X_d^V(g_{ab}v^b)-\frac{1}{2}X_a(g_{ij})v^iv^j-g_{ij}X_a^C(v^i)v^j -\frac{1}{2}X_a(g_{cd})v^cv^d -g_{cd}X_a^C(v^c)v^d\\
&=& v^{d} X_{d}\left(g_{a b}\right) v^{b}-v^{d} g_{a b} R_{d e}^{b} v^{e}+f^{d} g_{a b} \delta_{d}^{b}-\frac{1}{2} X_{a}\left(g_{c d}\right) v^{c} v^{d}  +g_{c d} R_{a e}^{c} v^{e} v^{d} 
\\
& =& g_{a b} f^{b}+\left(X_{d}\left(g_{a b}\right) v^{d} v^{b}-g_{a b} R_{d e}^{b} v^{d} v^e -\frac{1}{2} X_{a}\left(g_{c d}\right) v^{c}v^{d}  +g_{cd}R^c_{ae}v^ev^d\right) \\
& =& g_{a b} f^{b}+g_{a b} \Gamma_{ef}^{b} v^{e} v^{f}.
\end{eqnarray*}
We conclude that, indeed, $f^{b}=-\Gamma_{ef}^{b} v^{e} v^{f}$.

If, in the frame $\{X_\alpha\}$, the connection coefficients of the nonholonomic connection are given by  $\nabla^{nh}_{X_\alpha} X_\beta =
\tilde{\Gamma}_{\alpha \beta}^\gamma X_\gamma$, the corresponding quadratic spray can be expressed as $
\Gamma_{\nabla^{nh}}=v^\gamma X_\gamma^C-\tilde{\Gamma}_{\alpha \beta}^\gamma v^\alpha v^\beta X_\gamma^V$ (see e.g.\ \cite{reductionsodes}).
Specifically for the nonholonomic connection, we may compute:
\begin{itemize}
\item $\nabla^{nh}_{X_a}X_b=P(\Gamma_{ab}^\gamma X_\gamma)=\Gamma_{ab}^c X_c+0X_i$,

\item $\nabla_{X_a}^{nh} X_j=P\left(\Gamma_{a j}^\gamma X_\gamma\right)+\nabla_{X_a}^g X_j=\Gamma_{a j}^c X_c+\Gamma_{a j}^k X_k+\Gamma_{a j}^c X_c 
 =2 \Gamma_{a j}^c X_c+\Gamma_{a j}^k X_k $,

\item $\nabla_{X_i}^{nh} X_j=P\left(\Gamma_{ij}^\gamma X_\gamma\right)+\nabla_{X_i}^g X_j=\Gamma_{ij}^c X_{c} +\Gamma_{ij}^k X_k+\Gamma_{ij}^c X_c
=2 \Gamma_{ij}^c X_c+\Gamma_{ij}^k X_k$,
\end{itemize}
from which the  expression in the statement then follows.
\end{proof}

We now show that the second construction of a \sode\ extension is natural in the following way: it represents a generalization of the quadratic spray associated to the nonholonomic connection.
\begin{prop}  \label{prop:so}    
 For a purely kinetic Lagrangian $L(v)=\frac{1}{2}g(v,v)$,   the spray $\Gamma_{\nabla^{nh}}$ of the nonholonomic connection is the \sode\ extension $\Gamma^2$ of $\GLD$, i.e.\
$$
\Gamma_{\nabla^{nh}}=\Gamma_g-[\nabla(Q(\T))]^V.
$$
\end{prop}
\begin{proof}  When one restricts the expression of $\Gamma_{\nabla_{nh}}$ in the statement of Lemma~\ref{lem:VF} to $v^i=0$, one easily verifies that $\Gamma_{\nabla^{nh}}$ is indeed a \sode\ extension of $\GLD$. Besides, we can make the vector field $\Gamma_g$ appear in the expression of $\Gamma_{\nabla^{nh}}$, as follows:
\[   \Gamma_{\nabla^{nh}} = \Gamma_g-\Gamma^c_{\alpha j}v^\alpha v^j X_c^V+\Gamma^k_{\alpha a}v^\alpha v^a X_k^V.
\]

In case of the geodesic spray $\Gamma_g$, the coefficients of the dynamical covariant derivative are given  by $\Gamma^\beta_i=\Gamma^\beta_{i \alpha}v^\alpha$ and $\Gamma_a^\beta=\Gamma_{a \alpha}^\beta v^\alpha$. A more compact way of writing the spray of the nonholonomic connection is therefore 
$$
\Gamma_{\nabla^{nh}}= \Gamma_g-\Gamma^c_{j}v^j X_c^V+\Gamma^k_{ a} v^a X_k^V.
$$
For any quadratic spray (like $\Gamma_g$) is $\nabla\T=0$ (see e.g. \cite{derII}). For the projections $P$ and $Q$ (now interpreted as acting on $\mathcal{X}(\tau)$), we have
\begin{eqnarray*}
\nabla (Q(\T))&=& \nabla(Q\T) - Q(\nabla\T) = \nabla\left(v^i X_i\right)-Q\left(\Gamma^g\left(v^a\right) X_a+v^a \nabla X_a+\nabla\left(v^i X_i\right)\right) \\
& =&P\left(\nabla v^i X_i\right)-v^a Q\left(\nabla X_a\right) =v^{i}P\left(\nabla X_i\right)-v^a Q\left(\nabla X_a\right) 
=v^i P\left(\Gamma_i^\beta X_\beta\right)-v^a Q\left(\Gamma_a^\beta X_\beta\right) \\&=&v^i \Gamma_i^b X_b-v^a \Gamma_a^k X_k.
\end{eqnarray*}
When we take the vertical lift of this expression, it represents the difference between $\Gamma_{\nabla^{nh}}$ and  $\Gamma_g$, as required.
\end{proof}

\section{Geodesic extensions} \label{sec4}

In this section we consider only  purely kinetic nonholonomic  systems.  The \sode\ extension  $\Gamma_{\nabla^{nh}}$ (described in Proposition~\ref{prop:so}) is a quadratic spray, but it is not necessarily the geodesic spray of a (pseudo-)Riemannian metric $\hat g$. 

\begin{defi}  A geodesic extension of a purely kinetic nonholonomic system $(L,\D)$ is a Riemannian metric $\hat g$ whose geodesic spray $\Gamma_{\hat g}$ is a \sode\ extension of $\GLD$. 
\end{defi}

In that case, the nonholonomic trajectories could be interpreted as geodesics with initial velocities in $\mathcal{D}$. In the rest of the paper, we  look at this problem in more detail.

Consider a distribution $\D$ and a Riemannian metric $g$. We may use the metric $g$ to decompose
\[
TQ=\D \oplus \D^g.
\]
In what follows we will, again, make use of a frame $\{X_a,X_i\}$ that respects this decomposition and of its corresponding quasi-velocities $(v^a,v^i)$.

If $\hat g$ is a pseudo-Riemannian metric (i.e.\ non-degenerate, but not necessarily positive-definite), the integral curves of the geodesic spray $\Gamma_{\hat{g}}$ of $\hat{g} $ satisfy the equations
$$
\left\{\begin{array}{l}\dot{v}^a=-\hat{\Gamma}^a_{\alpha\beta}v^\alpha v^\beta=-\hat{\Gamma}_{b c}^a v^b v^c-\hat{\hat{\Gamma}}^a_jv^j , \\ 
\dot{v}^i=-\hat{\Gamma}^i_{\alpha\beta}v^\alpha v^\beta=-\hat{\Gamma}_{ab}^i v^a v^b -\hat{\hat{\Gamma}}^i_j v^j ,\end{array}\right.
$$
where $\hat{\Gamma}^\alpha_{\beta \gamma}$ are the connection coefficients of the Levi-Civita connection $\nabla^{\hat{g}}$ associated to the metric $\hat{g}$, in the frame $\{X_a,X_i\}$. {We have also used the notation $\hat{\hat{\Gamma}}^\alpha_jv^j=\hat{\Gamma}^\alpha_{jb}v^b+\hat{\Gamma}^\alpha_{ji}v^iv^j$ for all remaining terms that  include at least one factor $v^i$.} For our purpose, the terms with a factor  $v^j$ are of no special interest, since they vanish after restriction to $\mathcal{C}$. Recall that these equations are equivalent with the Euler-Lagrange equations of the Lagrangian $\hat{L}(v)=
\frac{1}{2}{\hat g}(v,v)= \frac{1}{2}\left(\hat{g}_{a b} v^a v^b+2 \hat{g}_{ai} v^{a} v^i+\hat{g}_{ij} v^i v^j\right)$.

When compared to the nonholonomic equations of motion of $(L,\D)$, the integral curves of any \sode\ extension of $\GLD$ should satisfy on $\C$:
$$
\left\{\begin{array}{l}\dot{v}^a=-\Gamma_{b c}^a v^b v^c, \\ \dot{v}^i=0.\end{array}\right.
$$ 
We conclude that the pseudo-Riemannian metric $\hat{g}$ we are looking for  must be such that its connection coefficients satisfy the  conditions 
$$
(*)\qquad\quad \left\{\begin{array}{l}\hat{\Gamma}_{ab}^c v^a v^b =\Gamma_{a b}^c v^a v^b,\\
\hat{\Gamma}_{ab}^i v^a v^b=0.\end{array}\right.
$$
{In the above indices $a$ and $b$ are summed over, and $i$ and $c$ are free indices. Also, the equations hold for any arbitrary choice of velocity coefficients $v^a$ and $v^b$.} Moreover, the factors $\hat{\Gamma}_{ab}^c$ are functions on $Q$ and do not depend on the quasi-velocities $v^a$. We may also write the  conditions $(*)$ in a version without the factors $v^a v^b$. However, we should keep  in mind  that (even though the Levi-Civita connection  is torsionless) the Christoffel symbols are not symmetric in the current anholonomic frame. For example, the last condition of $(*)$ is equivalent {to}
 \[
2\hat{\Gamma}_{ab}^i + R_{ab}^i= 0,
\]
and similar for the first. For compactness of our formulae, we prefer the first notation.

We can interpret the second condition in (*) also as follows: it shows that  $\Gamma_{\hat g}$ is  tangent to ${\mathcal C}$ because  then, on ${\mathcal C}=\{v^i=0\}$, \[
0=\Gamma_{\hat g}(v^i) \mid_{\mathcal C}= -{\hat\Gamma}^i_{\alpha\beta}v^\alpha v^\beta = - {\hat \Gamma}^i_{ab}v^a v^b.
\]

In what follows, we will make one further assumption: 
\begin{defi} A symmetric (0,2) tensor field $\hat g$ on $Q$ is a $\D$-preserving modification of a Riemannian metric $g$ on $Q$ if $\hat g$ and $g$ are related in such a way that
\[
{\hat g}|_{\D \times \D} = g|_{\D \times \D}.
\]
\end{defi}
In terms of the nonholonomic dynamics, this assumption has the interpretation that we are looking for an extension $\hat L$ of the mechanical Lagrangian $L$, where we do not wish to alter the constraint Lagrangian $L_c=L|_\C={\hat L}|_\C$.
In the current frame $\{X_a,X_i\}$, this assumption means that $g$ and $\hat g$ are of the form
$$
{g}=\left(\begin{array}{cc}
g_{a b} & 0 \\
0 & {g}_{i j}
\end{array}\right)
\quad \mbox{and}\quad
\hat{g}=\left(\begin{array}{cc}
g_{a b} & \hat{g}_{a i} \\ {(\hat{g}_{ai})^T} & \hat{g}_{i j}
\end{array}\right).
$$
In this, for example, $({\hat g}_{ij})$ is the matrix that corresponds to the restriction ${\hat g}|_{\D^g \times \D^g}$, i.e.\ ${\hat g}_{ij} = {\hat g}(X_i,X_j)$.

Recall that the principal restriction $g|_{\D \times \D}$ of a Riemannian metric remains positive-definite, because $g$ is as a whole. Among other, this means that the matrix $(g_{ab})$ with $g_{ab}=g(X_a,X_b)$ has an inverse $(g^{ab})$. Under the current assumptions,  this means that there exists also an orthogonal decomposition by $\hat g$,
\[
TQ = \D \oplus \D^{\hat g},
\]
where ${\rm Sec}(\D^{\hat g}) = {\rm span}\{ {\hat X}_i= X_i - {\hat g}_{bi} g^{ab} X_a \}$. In the next proposition, we will not assume from the outset that $\hat g$ is a Riemannian metric. The condition, in the Lemma below, that the restriction ${\hat g}|_{\D^{\hat g} \times \D^{\hat g}}$ is non-degenerate is equivalent with asking the matrix $({\hat g}({\hat X}_i,{\hat X}_j))=({\hat g}_{ij}-{\hat g}_{ai}{\hat g}^{ab}{\hat g}_{jb})$ to be invertible.

With these extra assumptions, we will show  that we can  re-write the conditions $(*)$ as algebraic and PDE conditions in the unknown functions ${\hat g}_{ai}$, which represents the restriction of $\hat g$ to $\D\times\D^g$. Besides, the Lemma below shows that the choice of ${\hat g}_{ij}$ is (almost) completely free. 

There is also a formulation of the conditions that solely makes use of  unknown functions $\theta_i={\hat g}_{ai}v^a$. In that case, also  the Lagrangian multipliers 
$\lambda_i=\Gamma_{(L,\mathcal{D})}(X_i^V(L))-X_i^C(L)$ (restricted to $\mathcal{C}$) make their appearance in the equations.  These are both components of a semi-basic 1-form on $TQ$. This is a  1-form that vanishes when applied to vertical lifts and which, for this reason, can be fixed by its action on complete lifts. Relying on the decomposition $TQ=\D\oplus\D^g$, we may consider the semi-basic 1-form defined by 
\[
\theta(X^C) = 0, \quad \forall X\in {\rm Sec}(\D), \qquad \theta(X^C)={\hat g}(\T,X) , \quad \forall X\in {\rm Sec}(\D^g).
\]
In the above, we have interpreted the Riemannian metric $\hat g$ again as the Hessian of its corresponding Lagrangian, i.e.\ as a (basic) (0,2) tensor field along $\tau$. We will be mainly interested in $\theta|_\C$, with components $\theta_i={\hat g}_{ai}v^a$, in the dual frame of $\{X_a^C,X_i^C,X_a^V,X_i^V\}$. 

Next, we consider the semi-basic 1-form, defined by 
\[
\lambda(X^C) = 0, \quad \forall X\in {\rm Sec}(\D), \qquad \lambda(X^C) =\Gamma_{(L,\mathcal{D})}(X^V(L))-X^C(L), \quad \forall X\in {\rm Sec}(\D^g), 
\]
and we consider its restriction to $\C$, $\lambda|_\C$. Again, its components are the restrictions of the multipliers
$\lambda_i=\Gamma_{(L,\mathcal{D})}(X_i^V(L))-X_i^C(L)$ to $\mathcal{C}$ we had encountered before.

\begin{lem}   \label{prop:result}
 Suppose that $g$ is a Riemannian metric and that $\hat g$ is a $\D$-preserving pseudo-Riemannian modification of $g$ whose restriction ${\hat g}|_{\D^{\hat g} \times \D^{\hat g}}$ is non-degenerate. Then, the equations (*)  are equivalent with the equations 
\begin{eqnarray*} (A) && \hat{g}_{bk} R^k_{ac}v^a v^b=0,  \\[2mm]
(B) && 
 \hat{g}_{di}{\Gamma}^d_{ab}v^a v^b=(g_{ki}\Gamma^k_{ab}+
X_a(\hat{g}_{bi}) -\hat{g}_{bk} R^k_{ai})v^a v^b,
\end{eqnarray*}
 in the unknown restriction  ${\hat g}|_{\D \times \D^g} = ({\hat g}_{ai})$. When written in terms of the 1-form $\theta = (\theta_i)$ the algebraic condition $(A)$ and the PDE condition $(B)$ become
\begin{eqnarray*} 
(A) && \theta_k R^k_{ac} v^a =0,\\
  (B) &&  \Gamma_{(L,\mathcal{D})}(\theta_i)+\theta_k R^k_{ia}v^a +\lambda_i=0.
	\end{eqnarray*}
\end{lem}

\begin{proof} 
In the frame $\{X_a,X_i\}$, we have $g_{ai}=0$.  Koszul's  formula (see Section~\ref{sec3}) for the metric $g$ leads to  
$$
\left\{\begin{array}{l}
2g_{c d} \Gamma_{a b}^dv^a v^b =(2X_a\left(g_{b c}\right)-X_c\left(g_{a b}\right)-2g_{d b} R_{a c}^d) v^av^b,\\
2g_{ki} \Gamma_{a b}^kv^a v^b =(-X_i\left(g_{ab}\right)-2g_{d b} R_{a i}^d) v^av^b.
\end{array}\right.
$$
Likewise, when applied on the metric $\hat{g}$ Koszul's formula gives
$$
\left\{\begin{array}{l}
2(\hat{g}_{c d} \hat{\Gamma}_{a b}^d+\hat{g}_{ck}\hat{\Gamma}^k_{ab})v^a v^b =(2X_a\left(\hat{g}_{b c}\right) -X_c\left(\hat{g}_{a b}\right)-2\hat{g}_{d b} R_{a c}^d 
-2\hat{g}_{bk}R^k_{ac}) v^av^b,\\
2(\hat{g}_{di} \hat{\Gamma}_{a b}^d+\hat{g}_{ki}\hat{\Gamma}^k_{ab})v^a v^b =(2X_a\left(\hat{g}_{b i}\right) -X_i\left(\hat{g}_{a b}\right)-2\hat{g}_{d b} R_{a i}^d 
-2\hat{g}_{bk}R^k_{ai}) v^av^b.
\end{array}\right.$$
Since $\hat{g}_{ab}=g_{ab}$,  we can recognize some terms of the first set in the second to obtain the following identities
\[
\begin{cases}
2({g}_{c d} \hat{\Gamma}_{a b}^d+\hat{g}_{ck}\hat{\Gamma}^k_{ab})v^a v^b =(2{g}_{c d} {\Gamma}_{a b}^d -2\hat{g}_{bk}R^k_{ac}) v^av^b,\\
2(\hat{g}_{di} \hat{\Gamma}_{a b}^d+\hat{g}_{ki}\hat{\Gamma}^k_{ab})v^a v^b =(2X_a\left(\hat{g}_{b i}\right)+2g_{ki} \Gamma_{a b}^k-2\hat{g}_{bk}R^k_{ai}) v^av^b.\end{cases}
\]

When conditions $(*)$ hold, the first identity above leads to condition $(A)$ and the second to condition $(B)$. 
We only need the extra assumptions in the statement to prove the converse. 

If we suppose that the conditions $(A)$ and $(B)$ hold, we get from the identities that
\[\begin{cases}
({g}_{c d} \hat{\Gamma}_{a b}^d+\hat{g}_{kc}\hat{\Gamma}^k_{ab})v^a v^b ={g}_{c d} {\Gamma}_{a b}^d v^av^b,\\
(\hat{g}_{di} \hat{\Gamma}_{a b}^d+\hat{g}_{ki}\hat{\Gamma}^k_{ab})v^a v^b =\hat{g}_{di}  {\Gamma}_{a b}^d  v^av^b.\end{cases}\nonumber
\]
We may rewrite the first as $\hat{\Gamma}_{a b}^dv^a v^b =  ({\Gamma}_{a b}^d - {g}^{c d}  \hat{g}_{ck}\hat{\Gamma}^k_{ab})v^av^b$. When plugged in the second, this gives $ (- \hat{g}_{di}{g}^{c d}  \hat{g}_{ck}+\hat{g}_{ki})\hat{\Gamma}^k_{ab}v^a v^b =0$. {Because we assume that $\hat{g}|_{\D^{\hat{g}}\times\D^{\hat{g}}}$ is non-degenerate,} we get indeed $\hat{\Gamma}^k_{ab}v^a v^b =0$ and therefore also $\hat{\Gamma}_{a b}^dv^a v^b =  {\Gamma}_{a b}^dv^av^b$.

For the statements involving $\theta_i$, it is clear that we may immediately rewrite $(A)$ in the desired form. For $(B)$, recall first that, on $\C$,
\begin{eqnarray*}
    \lambda_i&=& \Gamma_{(L,\mathcal{D})}(X_i^V(L))-X_i^C(L) = \Gamma_{(L,\mathcal{D})}(g_{ij}v^j)-\frac{1}{2}X_i(g_{ab})v^a v^b -g_{ab}X_i^C(v^a)v^b \\
    &=& g_{ij} \Gamma_{(L,\mathcal{D})}(v^j) -\frac{1}{2}X_i(g_{ab})v^a v^b -g_{ab}R^a_{ic}v^b = g_{ki}\Gamma^k_{ab}v^a v^b+g_{db}R^d_{ai}v^a v^b+g_{ab}R^a_{ic}v^c v^b \\&=& -g_{ki}\Gamma_L(v^k). 
\end{eqnarray*}
In {the second-to-last  equality}, we have again made use of Koszul's formula for the metric $g$.
Condition $(B)$ can then be expressed as
\begin{eqnarray*}
    0
&=& -\hat{g}_{di}\Gamma^d_{ab}v^a v^b+g_{ki}\Gamma^k_{ab}v^a v^b+
X_a(\hat{g}_{bi})v^a v^b-\hat{g}_{bk} R^k_{ai}v^a v^b 
\\&=&\hat{g}_{di}\Gamma_{(L,\mathcal{D})}(v^d) -g_{ki}\Gamma_L(v^k)+\Gamma_{(L,\mathcal{D})}(\hat{g}_{bi})v^b-\hat{g}_{bk}R^k_{ai}v^a v^b\\
&=& \Gamma_{(L,\mathcal{D})}(\hat{g}_{bi} v^b)-g_{ki}\Gamma_L(v^k)-\hat{g}_{bk}R^k_{ai}v^a v^b\\
&=& \Gamma_{(L,\mathcal{D})}(\hat{g}_{bi}v^b)+\lambda_i-\hat{g}_{bk}v^b R^k_{ai}v^a\\
&=&  \Gamma_{(L,\mathcal{D})}(\theta_i)+\lambda_i+\theta_k R^k_{ia}v^a.
\end{eqnarray*}
\end{proof}

{The conditions $(A)$ and $(B)$ are stated in terms of a local frame. The results in this paper are global only under the assumptions that such an adapted frame is globally defined.} The next proposition shows that the problem of finding a geodesic extension by means of a Riemannian metric is solved, once we have found a solution for the algebraic condition $(A)$ and the {PDE condition $(B)$}.

\begin{prop}   \label{prop:result2}
 Consider the equations $(A)$ and $(B)$ from the previous Lemma.
\begin{enumerate} \item If $\hat g$ is a Riemannian metric that is a $\D$-preserving modification of $g$ and a geodesic extension of $(L,\D)$, then its restriction ${\hat g}|_{\D\times\D^g}=({\hat g}_{ai})$ satisfies the algebraic conditions $(A)$ and the PDE conditions $(B)$.

\item If $\hat g$ is (0,2)-tensor field that is a $\D$-preserving modification of $g$ and whose restriction ${\hat g}|_{\D\times\D^g}=({\hat g}_{ai})$  satisfies the algebraic conditions $(A)$ and the PDE conditions $(B)$, then there exists a geodesic extension of $(L,\D)$ by a  Riemannian metric.
\end{enumerate}
\end{prop}

\begin{proof} Given what we have said in the proof of  Lemma~\ref{prop:result}, the first statement is obvious. To prove the second, we only need  to check that a restriction $\hat{g}|_{\D \times \D^g} = ({\hat g}_{ai})$ (solving $(A)$ and $(B)$) can always be completed to a full Riemannian metric $\hat g$, by an appropriate choice of the freedom in $\hat{g}|_{\D^g \times \D^g} = ({\hat g}_{ij})$. 

In this proof, we will use the shorthand notation $\hat{g}=\begin{pmatrix}
 {g}_{ab}& \hat{g}_{ai}\\\hat{g}_{ai}&\hat{g}_{ij}     
\end{pmatrix}=\begin{pmatrix}
    A&B\\B^T&D
\end{pmatrix}$.
Here, $A$ is a positive-definite $m\times m$-matrix, $B$ a $m\times k$-matrix and $D$ a $k\times k$-matrix. Below, we will  denote the spectrum of a matrix $A$ by $\sigma(A)$. Remark that if $\hat g$ represents a Riemannian metric (and therefore has $\det(D-B^T A^{-1}B)\neq 0$ in the current notations), also the condition that the restriction ${\hat g}|_{\D^{\hat g} \times \D^{\hat g}}$ should be non-degenerate in Lemma~\ref{prop:result} will be satisfied.

 If we consider the congruent  matrix $h=P{\hat g}P^T=\begin{pmatrix}
    A&0\\0&D-B^TA^{-1}B
\end{pmatrix}$, where $P= \begin{pmatrix}
    I_m&0\\B^TA^{-1}&I_k
\end{pmatrix}$
then,  by Sylvester's law of inertia, ${\hat g}$ and $h$ have the same signature. We may therefore conclude that, if we can make
$h$ positive definite, then so will be $\hat g$. The spectrum of $h$ is the union of the spectra of $A$
and $D - B^T A^{-1} B$, so our goal is to make $D - B^T A^{-1} B$ positive definite.  For that, we can choose $D=\alpha I_k$, with 
 $\alpha$ a yet to be determined function on $Q$. Since adding $\alpha I_k$ to a matrix shifts its spectrum
by $\alpha$, we can just take $\alpha$ in such a way that the spectrum of $-B^T A^{-1} B$ (which is a symmetric matrix with, therefore, a real spectrum)  is completely shifted into the
positives, that is, we can take any $\alpha$ such that $\alpha + {\rm min\,} \sigma(-B^T A^{-1} B) > 0$, or equivalently, $\alpha   > {\rm max\,}\sigma  (B^T A^{-1}B)$. Moreover, $\alpha$ can be chosen as a smooth (local) function on $Q$ since ${\rm max\,} \sigma (-B^T A^{-1} B)$ is a
continuous function. This ensures that $\hat g$ is smooth. 
\end{proof}

In Section~\ref{sec81} we will derive and integrate the conditions $(A)$ and $(B)$ explicitly for the example of a vertically rolling disk. We will also provide Riemannian metrics $\hat g$ which are geodesic extensions of the nonholonomic system.

The conditions $(A)$ and $(B)$, when written in terms of $\theta_i$ are familiar from a completely different context. If we set $\phi_i=-\theta_i$ (and forget that here $\theta_i$ is assumed to be linear in the $v^a$), the conditions
\[ 
\phi_k R^k_{ac} v^a =0 \quad\mbox{and}\quad \Gamma_{(L,\mathcal{D})}(\phi_i)+\phi_k R^k_{ia}v^a -\lambda_i=0
	\]
	also appear in \cite{anhol} (Corollary~1). The correspondence with our conditions  is only ostensibly, however, since there the conditions have a completely different interpretation and role in the theory: there, they are related to conditions for the nonholonomic problem to be consistent with a vakonomic problem. It would lead us too far  to explain here deeper the relation between this consistency problem and the geodesic extension problem we are now dealing with.\newline

 {We end this paragraph with a few comments on how to establish some information on the solvability of conditions $(A)$ and $(B)$. In practice, for a given nonholonomic system, the PDE condition $(B)$ can be hard to solve. The rank of the coefficient matrix in the algebraic condition $(A)$, on the other hand, determines the freedom of the candidate solutions $\theta_k$. In the next proposition, we  present a further necessary algebraic condition on the unknown $\theta_k$, in order for them to be a solution of both  $(A)$ and $(B)$.}

\begin{prop} \label{propsolv}
{If $\theta_k$ satisfies conditions $(A)$ and $(B)$, then it also satisfies  the  algebraic condition 
\[
\left[-S_i^k S^i_a+ \GLD(S^k_a)\right]\theta_k=\lambda_i S^i_a,
\]
where $S^i_a:=R^i_{ab}v^b$ and $S_i^k:=R_{ia}^k v^a$.}\end{prop}
\begin{proof}
{Suppose that  $\theta_k$ satisfies condition $(A)$, i.e. $\theta_i R^i_{ab}v^b=0$. When we apply the vector field $\GLD$ on both sides we get 
$$\GLD(\theta_i)R^i_{ab}v^b+\theta_i\GLD(R^i_{ab}v^b)=0. $$
If condition $(B)$ is satisfied as well, we may use  it to substitute the expression for $\GLD(\theta_i)$ in the first term. This gives 
$$\left(-\theta_k R^k_{ia}v^a-\lambda_i\right) R^i_{ab}v^b+\theta_i\GLD(R^i_{ab}v^b)=0, $$
and thus also $$\left[-S_i^k S^i_a+\GLD(S^k_a)\right]\theta_k=\lambda_i S^i_a.$$
}
\end{proof}

{Remark that we can repeat this technique over and over, to generate infinitely many algebraic conditions. They will only be new conditions, as long as they are independent from the previous ones. Because algebraic conditions are much more easy to work with than PDE conditions, in practical computations for a given system, one wishes to exploit first all these independent algebraic conditions. They allow one to narrow the possible candidates for $\theta_k$. Only in the last step one wants to check which ones of those $\theta_k$ solve the PDE condition $(B)$. We will show on a specific system how this can lead to an answer on the geodesic extension question, in the Section~\ref{carriagesec}.}

\section{Relation to \sode\ extensions} \label{sec5}

Let $(L=\frac{1}{2}g(v,v),\D)$ be a purely kinetic nonholonomic system. If we take a basis $\{X_a,X_i\}$ that respects the orthogonal decomposition by $g$ (i.e.\ with $\text{span}\{X_a\}=\mathcal{D}$ and $\text{span}\{X_i\}=\mathcal{D}^g$), we have already seen that the  nonholonomic vector field is
\[
\GLD=v^a X_a^C + f^a X_a^V = v^a X_a^C  -\Gamma_{bc}^av^bv^c X_a^V.
\]
Let   $\hat{g}$ be a Riemannian metric with $\hat{g}_{ab}=g_{ab}$. We have already mentioned that we may also use the metric  $\hat{g}$ to decompose $TQ$ into $TQ=\mathcal{D}\oplus {\D}^{\hat g}$. A basis  $\{\hat{X}_i\}$ that spans $ {\mathcal{D}}^{\hat g}$ is given by the vector fields $\hat{X}_i=X_i-{g}^{a b} \hat{g}_{b i} X_a$.
As in Section~\ref{sec3}, this decomposition will lead to a decomposition of the type 
\[
{\mathcal X}(TQ)={\mathcal X}^T(TQ)\oplus ({\mathcal X}^{\hat g}(\tau,\D))^V
\]
where  $({\mathcal X}^{\hat g}(\tau,\D))^V={\rm span}\{\hat{X}_i^V\}$.
Let's denote the projections that correspond to this  decomposition by   
$$ 
\hat{p}:\mathcal{X}(TQ)\rightarrow    {\mathcal X}^T(TQ) \quad\text{ and }\quad \hat{q}:\mathcal{X}(TQ) \rightarrow ({\mathcal X}^{\hat g}(\tau,\D))^V.
$$

The geodesic spray $\Gamma_{\hat g}$ (i.e.\ the  Lagrangian vector field $\Gamma_{\hat L}$ of $\hat{L}(v)=\frac{1}{2}\hat{g}(v,v)$) is given in terms of the basis $\{X_a,X_i\}$ by 
\[
 \Gamma_{\hat{g}}=v^\alpha X_\alpha^C+\hat{F}^a X_a^V+\hat{F}^i X_i^V = v^\alpha X_\alpha^C- {\hat{\Gamma}}^a_{\alpha\beta}v^\alpha v^\beta X_a^V-{\hat{\Gamma}}^i_{\alpha\beta}v^\alpha v^\beta X_i^V,
\]
where $\hat{\Gamma}^\gamma_{\alpha\beta}$ are the Christoffel symbols of $\hat{g}$, as before. 

Similarly as in Propositions~\ref{prop:1} and \ref{prop2}, we may construct from $\Gamma_{\hat{g}}$  two new \sode s. The first is 
\[
\Gamma_{\hat{g}}^1= \hat{p}(\Gamma_{\hat{g}}) =  v^\alpha X_\alpha^C+ ({\hat{F}}^a+g^{ab}\hat{g}_{bi}\hat{F}^i) X_a^V+0X_i^V.
\]

The second is (according to Proposition~\ref{prop:so}) the quadratic spray of the nonholonomic connection (now the one associated to $\hat g$ and $\D$). To derive its expression, we see that the orthogonal projection $\hat{Q}:\mathcal{X}(\tau)\rightarrow \mathcal{X}^{\hat g}(\tau, \mathcal{D})$ has the property
\[
\hat Q(\T)=\hat Q(v^a X_a+v^i X_i) = \hat Q((v^a+v^i g^{ab}\hat{g}_{bi}) X_a+v^i\hat{X}_i) = v^i {\hat X}_i .
\]
With that, we may define 
    \begin{eqnarray*}
        \Gamma_{\hat{g}}^2&=& \Gamma_{\hat{L}} -(\hat{\nabla}(\hat{Q}(\T)))^V\\
         &=& \Gamma_{\hat{g}}-(\Gamma_{\hat{g}}(v^i)\hat{X}_i+v^i\hat{\nabla}(\hat{X}_i))^V\\
        &=& v^\alpha X_\alpha^C+ ({\hat{F}}^a+g^{ab}\hat{g}_{bi}\hat{F}^i) X_a^V- v^i(\hat{\nabla}(\hat{X}_i))^V.
    \end{eqnarray*}

In the next proposition, we relate the condition $(A)$ to \sode\ extensions of $\GLD$ by one of these vector fields. 

\begin{prop} \label{propback1} Let $g$ and $\hat g$ be two Riemannian metrics which are a $\D$-preserving modification of each other.
 The following statements are equivalent:
 \begin{enumerate}
     \item[(i)] The condition $(A)$ is satisfied for $\hat g$,
     \item[(ii)] The vector field $\Gamma_{\hat{g}}^1$ is a \sode\  extension of $\Gamma_{(L_g,\mathcal{D})}$,
     \item[(iii)] The vector field $\Gamma_{\hat{g}}^2$ is a \sode\  extension of $\Gamma_{(L_g,\mathcal{D})}$,
     \item[(iv)] The nonholonomic vector field  $\Gamma_{(L_g,\mathcal{D})}$ satisfies on $\C$, 
		$$
		\Gamma_{(L,\mathcal{D})}(X^V(\hat{L}))-X^C(\hat{L}))=0, \qquad \forall X\in {\rm Sec}(\mathcal{D}).
		$$
 \end{enumerate}
\end{prop}

\begin{proof}  We have already remarked in the proof of Lemma~\ref{prop:result} that the Christoffel symbols ${\hat\Gamma}^\gamma_{\alpha\beta}$ of $\hat g$ and $\Gamma^\gamma_{\alpha\beta}$ of $g$ are always related as follows: 
$$
2({g}_{c d} \hat{\Gamma}_{a b}^d+\hat{g}_{kc}\hat{\Gamma}^k_{ab})v^a v^b =(2{g}_{c d} {\Gamma}_{a b}^d-\hat{g}_{bk}R^k_{ac}-\hat{g}_{ak}R^k_{bc}) v^av^b.
$$ 
The condition $(A)$ is satisfied if, and only if, we have
$$
2({g}_{c d} \hat{\Gamma}_{a b}^d+\hat{g}_{kc}\hat{\Gamma}^k_{ab})v^a v^b =(2{g}_{c d} {\Gamma}_{a b}^d) v^av^b.
$$
We may, equivalently, rewrite this as 
$$
g_{dc}\hat{F}^d|_{\C}+\hat{g}_{ck}\hat{F}^k|_{\C}=g_{dc}f^d.
$$

Since the matrix $g_{cd}$ is invertible we get 
$$
\hat{F}^d|_{\C}+g^{dc}\hat{g}_{ck}\hat{F}^k|_{\C}=f^d .
$$
Since the steps are reversible, this condition is in fact equivalent with the condition in (i). We continue the proof by showing that also all of the other statements are equivalent to this.

Indeed, given that $v^i=0$ on $\C$, it is easy to see from the expressions we had obtained before  that  
\[
\Gamma_{\hat{g}}^1|_\C= v^a X_a^C+ f^a X_a^V=\GLD
\]
if and only if $
\hat{F}^d|_{\C}+g^{dc}\hat{g}_{ck}\hat{F}^k|_{\C}=f^d
$, and similarly for $\Gamma_{\hat{g}}^2$.

The statement $(iv)$ says that the coefficients $f^a$ of the nonholonomic vector field $\GLD$ satisfy, on $\mathcal{C}$, 
\begin{eqnarray*}
  0  &= &\Gamma_{(L_g,\mathcal{D})}(X_a^V(\hat{L}))-X_a^C(\hat{L})\\
    &=& v^b X_b^C X_a^C(\hat{L})+f^b X_b^V X_a^V(\hat{L})-X_a^C(\hat{L})\\
    &=& v^b X_b^C X_a^C(\hat{L})+f^b {g}_{ab}-X_a^C(\hat{L}).
\end{eqnarray*}

Now, recall that the Lagrangian vector field $\Gamma_{\hat{L}} = \Gamma_{\hat{g}} $ always satisfies
$$
\Gamma_{\hat{L}}(X^V(\hat{L}))-X^C(\hat{L})=0,\quad\forall X\in\mathcal{X}(Q).
$$
In particular, we have (after restriction to $\C$)
\begin{eqnarray*}
  0 & =&\Gamma_{\hat{L}}(X_a^V(\hat{L}))-X_a^C(\hat{L}) = v^bX_b^CX_a^V(\hat{L})+\hat{F}^b|_\mathcal{C}X_b^V X_a^V(\hat{L})+F^i|_{\mathcal{C}}X_i^V X_a^V(\hat{L})-X_a^C(\hat{L}) \\
    &=&v^bX_b^CX_a^V(\hat{L})+\hat{F}^b|_\mathcal{C}{g}_{ab}+F^i|_{\mathcal{C}}\hat{g}_{ai}-X_a^C(\hat{L}). 
\end{eqnarray*}

If we subtract this from the previous, the statement in $(iv)$ is indeed equivalent with $\hat{F}^b|_\mathcal{C}+F^i|_{\mathcal{C}}\hat{g}_{ai}g^{ab}=f^b$.
\end{proof}

We have already shown in Proposition~\ref{prop:result2}  that $\Gamma_{\hat{g}}$ is a \sode\  extension of $\Gamma_{(L_g,\mathcal{D})}$ if and only if $(A)$ and $(B)$ are satisfied. From the previous we know that, if only $(A)$ is satisfied, it is not enough to guarantee that $\Gamma_{\hat g}$ is a \sode\ extension: Proposition~\ref{propback1} says that in that case only $\Gamma_{\hat g}^1$ or $\Gamma_{\hat g}^2$ is one. Therefore, we need to add the condition that also the restriction
$\hat{\nabla}(\hat{Q}(\T))|_{\mathcal{C}}$ vanishes. The proof of the next proposition is now  obvious.

\begin{prop} \label{propback2}
Let $g$ and $\hat g$ be two Riemannian metrics which are a $\D$-preserving modification of each other.    The following statements are equivalent:
    \begin{enumerate}
\item[(i)] $\Gamma_{\hat{g}}$ is a SODE extension of $\Gamma_{(L,\mathcal{D})}$,
        \item[(ii)] The conditions $(A)$ and  $(B)$ are satisfied,
\item[(ii)] (A) is satisfied and $\hat{\nabla}(\hat{Q}(\T))|_{\mathcal{C}}=0$,
\item[(iii)] $\Gamma_{\hat{g}}|_\mathcal{C}=\Gamma_{\hat{g}}^1|_\mathcal{C} = \Gamma_{\hat{g}}^2|_\mathcal{C}$.
    \end{enumerate}
\end{prop}

\section{Infinitesimal Gauss-type conditions}\label{sec6}

Throughout this section, given $v_q\in T_qQ$, we will often make use of the canonical identification of $T_{v_q}(T_q Q)$ as the tangent space $T_q Q$: 
$$
(w_{q})_{{v_q}}^V=w^\alpha X_\alpha^V(v_q) \,\,\longleftrightarrow\,\, w_q=w^\alpha X_\alpha(q).
$$

Let $\phi_t^{nh}: \D\to\D$ be the flow of the nonholonomic vector field $\GLD$. The nonholonomic exponential map is defined in \cite{firstnonholexp} as \[
\exp^{nh}_{q}(v_q) = \tau(\phi_1^{nh}(v_q)),
\]
where $\tau:TQ\to Q$ is the tangent bundle projection.
Among other properties, it has been shown in \cite{Simoes} specifically for a purely kinetic nonholonomic system $(L(v)=\frac{1}{2}g(v,v),\D)$ that for each fixed $q\in Q$, there exists a submanifold $M_q$ of $Q$ (with $q\in M_q$) such that the nonholonomic exponential map restricts to a diffeomorphism 
$\exp^{nh}_q: U_0 \subset \D_q \to  M_q$, where $U_0$ is a starshaped
open subset of {$\D_q$} around $0_q \in U_0$.  Moreover,   $\exp^{nh}_q(0_q) = q$ and  the  morphism $T_{0_q} (\exp^{nh}_q)
: \C_q \to T_qQ$ (under the identification between $\C_q$ and $T_{0_q}U_0$) is the canonical inclusion. Important for us, is that it has also  been shown in \cite{Simoes} that, for every $v_q \in U_0$,
\[
\exp^{nh}_q (tv_q) = c_{v_q}(t),
\]
with $c_{v_q}: [0, 1] \to M_q$ the nonholonomic trajectory satisfying $c_{v_q}(0) = q$ and ${\dot c}_{v_q}
(0) = v_q$.

In what follows we will use a description of the nonholonomic exponential in quasi-velocities. Given a frame $\{X_\alpha\} = \{X_a,X_i\}$, adapted to $TQ=\D \oplus \D^g$, as before, we set  in local coordinates
$$
\exp^{nh}_q:  v_q=v^a X_a(q)\mapsto \tilde{q}^\alpha =\epsilon^\alpha(v).
$$
We will also  use $((E^\alpha_\beta)(v))$ for the matrix representation of the tangent map of $\exp^{nh}_q$ in $v_q$ in quasi-velocities: For $u_q\in \D_q$, we set 
\[
    T_{v_q}(\exp^{nh}_q)(u_q)  = u^a \left(E^b_a(v)X_b(\tilde{q})+E^i_a(v) X_i(\tilde{q})\right) \in T_{\tilde q}Q.
\]

\begin{lem} \label{hulplemma} At $0_q\in T_qQ$, we have 
$E_a^b(0)=\delta_a^b$, $E_a^i(0)=0$,
$$
\frac{\partial E^d_a}{\partial v^b}(0)=\frac{1}{2}R^d_{ab}(q)-\frac{1}{2}(\Gamma_{ab}^d(q)+\Gamma^d_{ba}(q)) \quad\text{and}\quad \frac{\partial E^i_a}{\partial v^b}(0)=\frac{1}{2}R^i_{ab}(q).
$$
\end{lem}

\begin{proof} The first two properties follow from the fact that, when $v_q=0_q$,   $T_{0_q}(\exp^{nh}_q)$ can be regarded as the inclusion.
We will compute the other expressions in two steps. 

Let $v_q = v^aX_a(q)\in T_qQ$. Consider a nonholonomic trajectory $c_{v_q}(t)$ and its derivative $\dot{c}_{v_q}(t)=V^a(t) X_a\left(c_{v_q}(t)\right) \in \mathcal{D}_{c_{v_q}(t)}$. The functions $V^a(t)$ then satisfy $\dot{V}^a(t)=-\Gamma_{b c}^a\left(c_{v_q}(t)\right) V^b(t) V^c(t)$.
Since in the purely kinetic case {$c_{v_q}(t)=\exp^{nh} \left(t v_q\right)$}, we may also write this as 
\[
\dot{c}_{v_q}(t)=E_a^\alpha(tv) v^a X_\alpha\left(c_{v_q}(t)\right).
\]
We get from this  
$$
E_{a}^d(t v_q) v^a=V^d(t) \quad \text{and} \quad E_a^i(t v_q) v^a=0.
$$
By applying $\frac{d}{dt}|_{t=0}$ on both sides of both equations we get 
$$
\frac{\partial E_a^d}{\partial v^\mu}(0) v^\mu v^a= \dot{V}^d(0)=-\Gamma^d_{ef}(q) v^e v^f
$$ and after symmetrization
$$ \frac{\partial E_a^d}{\partial v^c}(0)+\frac{\partial E_c^d}{\partial v^a}(0) =-(\Gamma_{ac}^d(q)+\Gamma_{ca}^d(q)).
$$
Likewise, it follows from 
$
{\partial E_a^i}/{\partial v^\mu}(0) v^\mu v^a = {\partial E_a^i}/{\partial v^b}(0) v^b v^a=0
$ that 
${\partial E_a^i}/{\partial v^b}(0)+{\partial E_b^i}/{\partial v^a}(0)=0$.

Now we compute the expression for the differences $ {\partial E_a^d}/{\partial v^c}(0)-{\partial E_c^d}/{\partial v^a}(0) $ and ${\partial E_a^i}/{\partial v^b}(0)-{\partial E_b^i}/{\partial v^a}(0)$. 
We know that in natural coordinates and in quasi-velocities, we get (respectively) 
\[
    T_{v_q}(\exp^{nh}_q)(u_q) = u^a\frac{\partial \epsilon^\alpha}{\partial v^a}(v) \frac{\partial}{\partial q^\alpha}|_{\tilde{q}} = E_a^\beta(v) u^a X_\beta(\tilde{q}).
		\]
Since $\{{\partial}/{\partial q^\alpha}|_{\tilde q}\}$ is a basis for $T_{\tilde q} Q$ there  exists an invertible matrix  $(A_\beta^\alpha(\tilde{q}))$ of functions such that
$X_\beta(\tilde{q})=A_\beta^\alpha(\tilde{q}) {\partial}/{\partial q^\alpha}|_{\tilde{q}}$. Therefore 
$ {\partial \epsilon^\alpha}/{\partial v^a}(v)=E_a^\beta(v) A_\beta^\alpha(\epsilon(v))$, and  in particular ${\partial \epsilon^\alpha}/{\partial v^a}(0)=E_a^\beta(0) A_\beta^\alpha(q)=A_a^\alpha(q)$.
When we take a $\partial/\partial v^b$-derivative on both sides of the first relation, we get
\[
    \frac{\partial^2 \epsilon^\alpha}{\partial v^a\partial v^b}(v)=\frac{\partial E_a^\beta}{\partial v^b}(v) A_\beta^\alpha(\epsilon(v))+E_a^\beta(v) \frac{\partial A_\beta^\alpha}{\partial q^\mu}(\epsilon(v))\frac{\partial \epsilon^\mu}{\partial v^b}(v),
	\]
and at  $v=0$:
\[
\frac{\partial^2 \epsilon^\alpha}{\partial v^a\partial v^b}(0)=\frac{\partial E_a^\beta}{\partial v^b}(0) A_\beta^\alpha(q)+\delta_a^\beta \frac{\partial A_\beta^\alpha}{\partial q^\mu}(q)
A^\mu_b(q).
\]
If we change in this last equality  $a$ and $b$, and subtract we get the identity
$$
 0=\left(\frac{\partial E_a^\beta}{\partial v^b}(0)-\frac{\partial E_b^\beta}{\partial v^a}(0)\right) A_\beta^\alpha(q)+\left( \frac{\partial A_a^\alpha}{\partial q^\mu}(q)A^\mu_b(q)- \frac{\partial A_b^\alpha}{\partial q^\mu}(q)A^\mu_a(q) \right).
$$
The last term is, in fact, related to the Lie bracket
$ R^\beta_{ba} X_\beta= [X_b,X_a]=(A^\mu_b{\partial A_a^\alpha}/{\partial q^\mu}- A^\mu_a{\partial A_b^\alpha}/{\partial q^\mu}) {\partial}/{\partial q^\alpha}$. 
The difference now becomes
$$ 
0=\left[\frac{\partial E_a^\beta}{\partial v^b}(0)-\frac{\partial E_b^\beta}{\partial v^a}(0) -R_{ab}^\beta(q)\right]A^\alpha_{\beta}(q).
$$
For $\beta=d$ and $\beta=i$, we get, respectively,
 $$
\frac{\partial E_a^d}{\partial v^b}(0)-\frac{\partial E_b^d}{\partial v^a}(0)=R_{ab}^d(q) ,
\qquad \frac{\partial E_a^i}{\partial v^b}(0)-\frac{\partial E_b^i}{\partial v^a}(0) =R_{ab}^i(q).
$$
We conclude that 
\begin{equation}
    \begin{cases}
\displaystyle      \frac{\partial E_a^d}{\partial v^b}(0)-\frac{\partial E_b^d}{\partial v^a}(0)=R_{ab}^d(q),\\[2mm]
   \displaystyle   \frac{\partial E_a^d}{\partial v^c}(0)+\frac{\partial E_c^d}{\partial v^a}(0) =-(\Gamma_{ac}^d(q)+\Gamma_{ca}^d(q)),
    \end{cases}\nonumber  \begin{cases}
\displaystyle        \frac{\partial E_a^i}{\partial v^b}(0)-\frac{\partial E_b^i}{\partial v^a}(0) =R_{ab}^i(q), \\[2mm]
\displaystyle        \frac{\partial E_a^i}{\partial v^b}(0)+\frac{\partial E_b^i}{\partial v^a}(0)=0 .
    \end{cases}
\end{equation}
This leads to the expressions in the statement of the Lemma.
\end{proof}

We can now relate condition $(A)$ to the Gauss Lemma of Riemannian geometry.
Let $h$ be a Riemannian metric on a manifold $Q$ with exponential map $\exp^h_q$. Suppose that $v_q\in T_qQ$ is such that $\exp_q^h(v_q)$ is defined. 
The Gauss Lemma says, in essence, that  $\exp_q^h:T_qQ \to Q$ is  a radial isometry (under appropriate identification). This means that
${\mathcal H}^0=(\exp^{h}_q)^*h$ satisfies the `Gauss condition':
$$
\mathcal{H}^0_{v_q}(v_q,u_q)=\mathcal{H}^0_{0_q}(v_q,u_q),
$$ 
for all  $u_q\in T_{v}T_{q}Q\cong T_{q}Q$.
The expression for $\mathcal{H}^0$ can equivalently be written as
\[
h_{\exp^h_{q}(v_q)}(T_{v_q}\exp^h_{q}(v_q),T_{v_q}\exp^h_{q}(u_q))=h_q( v_q,u_q).
\]
 We refer to e.g.\ \cite{do} for the details. 

We do not need the Gauss Lemma as such. In \cite{Simoes} it is shown that a certain Gauss-type condition is necessary and sufficient for the existence of a Riemannian metric $g^{nh}_q$ on the submanifold $M_q$ with the property that the radial kinetic nonholonomic trajectories departing from $q$ are minimizing geodesics of $g^{nh}_q$. This is also not the property that we will be interested in here, but our result has a similar flavour. In the next proposition we give a sufficient condition under which the condition $(A)$ (that appears in Lemma~\ref{prop:result} and Proposition~\ref{prop:result2}) always holds. 

\begin{prop} \label{prop:gauss} Let $\exp^{nh}_q$ be the nonholonomic exponential map of a purely kinetic $(L,\D)$ and let $\hat{g}$ be a Riemannian metric on $Q$ which is a $\D$-preserving modification of $g$. If $\mathcal{G}^0=(\exp^{nh}_q)^* \hat{g}$ satisfies  the following Gauss-type condition 
$$
\mathcal{G}^0_{v_q}(v_q,u_q)=\mathcal{G}^0_{0_q}(v_q,u_q),\quad\quad\forall v_q,u_q\in\mathcal{D}_q,
$$
or, equivalently, 
$$
\hat{g}_{\exp^{nh}_q(v_q)}\left(T_{v_q}\left(\exp _q^{nh}\right)\left(v_{q}\right), T_{v_q}\left(\exp _q^{nh}\right)\left(u_{q}\right)\right)=\hat{g}_q(v_{q},u_{q}),
$$
 then condition $(A)$  is satisfied for $\hat{g}$. 
\end{prop}

\begin{proof} 
The equivalent condition is the first written in full, after having taken into account that the tangent map of the exponential map {at} zero is the inclusion (with the identification between $T_q Q$ and $T_{v_q}(T_q Q)$).

Given that  $v_q,u_q\in\mathcal{D}_q$, we can write them as  $v_q=v^a X_a(q)$ and $w_q=u^a X_a(q)$. With this the Gauss condition becomes
$$
\hat{g}_{\alpha\beta}(\exp^{nh}_q(v_q)){E_a^\alpha(v)} v^a E_b^\beta(v)=g_{a b}(q) v^a.
$$

Now, we substitute $v_q$ by $t v_q$ and cancel a factor $t$ on both sides, 
$$
 \hat{g}_{{\alpha\beta}}(c_{v_q}(t))E_a^\alpha(tv)v^a E_b^\beta(t v)=g_{a b} v^a,
$$
where we have also used $\exp^{nh}_q(tv_q)=c_{v_q}(t)$.
If we apply $\left.\frac{d}{d t}\right|_{t=0}$ on both sides, we obtain 
\begin{eqnarray*}
    \frac{\partial \hat{g}_{\alpha\beta}}{\partial q^\mu}(q)(\dot{c}_{v_q})^\mu(0) {E_a^\alpha(0)} v^a {E_b^\beta(0)} +\hat{g}_{{\alpha\beta}}(q) \frac{\partial E_a^\alpha}{\partial v^\mu}(0) v^\mu v^a {E_b^\beta(0)} +\hat{g}_{{\alpha\beta}}(q)  {E_a^\alpha(0)} v^a \frac{\partial E_b^\beta}{\partial v^\mu}(0) v^\mu=0.
\end{eqnarray*}
In the first term, we may rewrite $(\dot{c}_{v_q})^\mu(0){\partial}/{\partial q^\mu}|_q$ as $v^f X_f(q)+v^i X_i(q)=v^f X_f (q)$, since the vector   $v_q = {\dot c}_{v_q}(0)\in\mathcal{D}_q$ has components $v^i=0$. If we also insert the expressions that we had found in Lemma~\ref{hulplemma}, we get
\begin{eqnarray*}
    0 
     &=& v^f X_f(\hat{g}_{ab})v^a+\hat{g}_{{d b}}\left( -(\Gamma^d_{af}+\Gamma_{fa}^d)\right) v^f v^a+ \frac{1}{2}\hat{g}_{{ad}}\left( -(\Gamma^d_{bf}+\Gamma^d_{jb})+R^d_{bf}\right) v^f v^a+
   \frac{1}{2}\hat{g}_{{ai}}R^i_{bf}v^f v^a
	\\
   &=&  \left( X_f(\hat{g}_{ab})-\frac{1}{2}\hat{g}_{{d b}}\big(2\Gamma^d_{af}+0\big)-\frac{1}{2}\hat{g}_{{ad}}\big(2\Gamma^d_{fb}\big)+
   \frac{1}{2}(\frac{1}{2}(\hat{g}_{ai}R^i_{bf}+\hat{g}_{fi}R^i_{ba}))\right)v^av^f\\
   &=& \left(X_f(\hat{g}_{ab})-\hat{g}_{db}\Gamma^d_{af}-\hat{g}_{ad}\Gamma^d_{fb}+\frac{1}{4}(\hat{g}_{ai}R^i_{bf}+\hat{g}_{fi}R^i_{ba})\right)v^av^f.
\end{eqnarray*}
Recall first that for the metric $g$, $\nabla^g_{X_f}g=0$. In the current frame this gives 
 $$
\left(X_f({g}_{ab})-\hat{g}_{db}\Gamma^d_{af}-{g}_{ad}\Gamma^d_{fb}\right)v^a v^f =0.
$$
Since we assume that $\hat{g}_{ab}=g_{ab}$, we may use this to see  the first terms in the previous expression  vanish. What remains is then indeed condition $(A)$ from Lemma~\ref{prop:result}, 
$$
\left(\hat{g}_{ai}R^i_{bf}+\hat{g}_{fi}R^i_{ba}\right)v^av^f=0,
$$
or, equivalently,
$
\hat{g}_{ai}R^i_{bf}v^av^f=0$.
\end{proof}

 We end this paragraph with an equivalent characterization of the   condition $(B)$.
It will be most convenient to come back to the version of $(B)$, in the form where we regard it as a PDE in  the unknowns $\theta_i={\hat g}_{ai}v^a$ (i.e.\ the restriction of some semi-basic 1-form $\theta$ to $\C$). 
We now show that  $(B)$ can be interpreted as the infinitesimal version of a geometric property  along $c_{v_q}(t)$. To do so, we introduce yet {another} \sode\ extension $\Gamma_{\overline\nabla}$ of $\GLD$, which is the quadratic spray  of a linear connection $\overline\nabla$ on $Q$, different from the nonholonomic connection $\nabla^{nh}$. The construction of the connection $\overline\nabla$ is inspired by a result that can be found in \cite{PF}.

We have already introduced the two projection operators $P:TQ\rightarrow \mathcal{D}$ and $Q: TQ\rightarrow \mathcal{D}^g$, related to $g$. Consider the Levi-Civita connection $\nabla^g$ of $g$. We may use it to define the operator $\nabla^P$ as
$$
\nabla^P: {\rm Sec}(\mathcal{D})\times {\rm Sec}(\mathcal{D})\rightarrow\mathcal{X}(Q),\quad (X,Y)\mapsto \nabla^P_X Y= {P}(\nabla^g_X Y).
$$
Proposition~2.1 in \cite{PF} states that we may extend this  to a covariant derivative
$$
\overline{\nabla}^P:\mathcal{X}(Q)\times {\rm Sec}(\mathcal{D}) \rightarrow\mathcal{X}(Q),\quad (X,Y)\mapsto \nabla^P_{P(X)}Y+P([Q(X),Y]).
$$
We can do the same for the projection operator $Q$. Next, from Theorem~2.2 in \cite{PF}, we obtain that, for $X,Y\in\mathcal{X}(Q)$, the operator
$$
\overline\nabla:\mathcal{X}(Q)\times \mathcal{X}(Q)\rightarrow \mathcal{X}(Q),\quad (X,Y)\mapsto \overline{\nabla}^P_X({P}(Y))+\overline{\nabla}^Q_X({Q}(Y))
$$
is a linear connection on $Q$. In the current frame, the coefficients of this linear connection  are
$$
{\overline\nabla}_{X_a}X_b=\Gamma_{ab}^c X_c,\quad {\overline\nabla}_{X_a}X_i=R^k_{ai} X_k,\quad {\overline\nabla}_{X_i} X_a=R^c_{ia} X_c,\quad {\overline\nabla}_{X_i}X_j=\Gamma_{ij}^k X_k,
$$ 
and its corresponding spray is 
$$
\Gamma_{\overline\nabla}=v^\alpha X_\alpha^C+(\Gamma_{ab}^c v^a v^b +R^c_{ia}v^i v^a )X_c^V+ (R^k_{ai} v^a v^i +\Gamma_{ij}^k v^i v^j) X_k^V.
$$
By restricting it to $v^i=0$ we can easily see that it is, indeed, a \sode\ extension of $\GLD$.

Let $c_{v_q}(t)$ be a nonholonomic trajectory, with
${\dot c}_{v_q}(t) = V^a(t)X_a(c_{v_q}(t)) \in \D_{c_{v_q}(t)}$.
Since $\Gamma_{\overline\nabla}$ is a \sode\ extension of $\GLD$, we may  interpret $c_{v_q}(t)$ as one of its base integral curves. We may then use the linear connection ${\overline\nabla}$ to parallel transport a vector $w_q\in T_qQ$ along $c_{v_q}(t)$. In this way, we obtain a vector field $W(t)$ along $c_{v_q}(t)$ whose components along $\{X_a,X_i\}$ are solutions of the initial value problem  
$$
\dot{W}^\alpha +{\overline\Gamma}^\alpha_{\gamma\beta}V^\gamma W^\beta=0, \qquad W^\alpha(0) = w^\alpha.
$$
For our specific case, we get
\begin{equation*}
    \begin{cases}
       \dot{W}^c=-{\Gamma}_{ab}^cV^a W^b,\\
       \dot{W}^i=-R_{aj}^iV^a W^j.
    \end{cases}\nonumber
\end{equation*} 
The above differential equations are separated. If we take the initial value $w_q\in \D_q^g$, we have $w^a=0$ and therefore is the solution of the first simply $W^a(t)=0$. This means that $W(t) \in \D_{c_{v_q}(t)}^g$.

For a 1-form $h(t)$ along the curve $c_{v_q}(t)$, its covariant derivative ${\overline\nabla}_{{\dot c}_{v_q}(t)} h(t)$ is again a 1-form along $c_{v_q}(t)$. It can be defined below by means of its action on vector fields $X(t)$ along  $c_{v_q}(t)$:
\[
\Big(  {\overline\nabla}_{{\dot c}_{v_q}(t)} h(t) \Big) \Big(X(t)\Big) = \frac{d}{dt} \Big( \Big( h(t) \Big) \Big(X(t)\Big)\Big)  -   \Big( h(t) \Big) \Big( {\overline\nabla}_{{\dot c}_{v_q}(t)}  X(t) \Big).
\]

Consider $\{X^i, X^a\}$ the basis of 1-forms, dual to the basis $\{X_i, X_a\}$ of vector fields. Then ${\rm span}\{X^i\}$ is the annihilator of $\mathcal{D}$, say $\mathcal{D}^o$.
Since $\dot{c}_{v_q}(t)$ has only components in $\mathcal{D}_{c_{v_q}(t)}$, we can compute that, along $c_{v_q}(t)$, the covariant derivative of $h(t)=h_i(t)X^i(c_{v_q}(t))+h_b(t)X^b(c_{v_q}(t))$ is given by
\[
{\overline\nabla}_{\dot{c}_{v_q}(t)}h = \left(\dot{h}_c- {\Gamma}^b_{a c}V^a h_{b} \right)X^c + 
\left(\dot{h}_i- {R}^k_{a i}V^a  h_{k}\right)X^i.
\]
Now, consider the  1-form $\lambda(t) = \lambda({\dot c}_{v_q}(t))$ along  $c_{v_q}(t)$, i.e.\ the restriction of the semi-basic 1-form $\lambda$ that represented the Lagrangian multipliers (in Section~\ref{sec4}). For this, we can define a unique 1-form  $h(t)$ along $c_{v_q}(t)$ that satisfies 
   \[      {\overline\nabla}_{\dot{c}_{v_q}}h(t)=\lambda(t),\qquad h(0)=0_q\in T^*_qQ.
\]
This $h(t)$ is the unique solution of the initial value problem 
\[
    \begin{cases}
     \dot{h}_c(t)=\Gamma^b_{ac}(c_{v_q}(t))V^a(t) h_b(t),\\ \dot{h}_i(t)=R^j_{ai}(c_{v_q}(t))V^a(t) h_j(t)+\lambda_i(t),
		\end{cases}\qquad
		\begin{cases}
     h_c(0)= 0,\\
     h_i(0)=  0.
    \end{cases}
		\]
From the first relation, we may conclude that $h_b(t)=0$.

We have now gathered all the tools to define an equivalent characterization of the PDE condition $(B)$ from Proposition~\ref{prop:result2}.

\begin{prop} \label{prop:PT} Consider $q\in Q, v_q\in\mathcal{D}_q$ and the 1-form $h(t)$, as defined above. A $\D$-preserving modification  $\hat{g}$ of $g$ satisfies $(B)$ if and only if  
$$ 
\hat{g}_{c_{v_q}(t)} \left( {\dot c}_{v_q}(t),W(t)\right)+ \Big(h(t)\Big)\Big(W(t)\Big)=\hat{g}_q(v_q,w_q), \qquad \forall  w_q\in\mathcal{D}^g_q,
$$
   where 
	$W(t)$ is the parallel transport of $w_q\in\mathcal{D}^g_q$ with respect to ${\overline\nabla}$.
	\end{prop}

Apart from the term in $h(t)$ this condition is again reminiscent to a Gauss-type condition.
If we multiply the relation in the proposition on both sides with an extra factor $t$, we can  express it equivalently  in terms of the exponential map $\exp^{nh}_q(tv_q)=c_{v_q}(t)$, as follows:
$$ 
\hat{g}_{\exp_q^{nh}(tv_q)}\left(T_{tv_q}(\exp^{nh}_q)(tv_q),W(t)\right)+t \Big(h(t)\Big)\Big(W(t)\Big)=\hat{g}_q(tv_q,w_q).
$$

\begin{proof} 
From Lemma~\ref{prop:result} we know that the expression $(B)$ can be written as 
\[
\Gamma_{(L,\mathcal{D})}(\theta_i)+\lambda_i+\theta_k R^k_{ib}v^b=0,
\]
or, after expanding and multiplying both sides with $w^i$: 
\[
X_f(\hat{g}_{ai})v^f v^a w^i+  {\hat g}_{ai}f^aw^i  +\lambda_i w^i +  \hat{g}_{ka}v^a R^k_{ib}v^b w^i=0.
\]
This expression is valid for arbitrary choices of $q,v^a$ and $w^k$. If we replace them, specifically, by $c_{v_q}(t)$, $V^a(t)$ and $W^k(t)$, respectively, we get, for all $t$:
\[
X_f(\hat{g}_{ai})V^f(t) V^a(t) W^i(t)+  {\hat g}_{ai}{\dot V}^a(t) W^i(t)  +\lambda_i(t) W^i(t) +  \hat{g}_{ak}V^a(t) R^k_{ib}V^b(t) W^i(t)=0
\]
or, 
\[
\frac{d}{dt}\left( \hat{g}_{ai}  V^a(t) \right) W^i(t)+  \left( \dot{h}_i(t)-R^j_{ai}(c_{v_q}(t))V^a(t) h_j(t)   \right) W^i(t) +  \hat{g}_{ak}V^a(t) {\dot W}^k(t)=0.
\]
We may simplify this to
\[
\frac{d}{dt}\left( \hat{g}_{ai}  V^a(t)  W^i(t)\right)+  \left( \dot{h}_i(t) W^i(t)  + {\dot W}^j(t) h_j(t) \right)   =0,
\]
from which
\[
\frac{d}{dt}\left( \hat{g}_{ai}  V^a(t)  W^i(t)+ {h}_i(t) W^i(t) \right)   =0.
\]
This means that 
\[
\hat{g}_{ai}  V^a(t)  W^i(t)+ {h}_i(t) W^i(t)=\hat{g}_{ai}  v^a  w^i + {h}_i(0) w^i.
\]
Since $h_i(0)=0$, we get the result.

The other direction follows if we start from the expression above, take a derivative by $t$ and set $t=0$, afterward.
\end{proof}

From Propositions~\ref{prop:result2}, \ref{prop:gauss} and \ref{prop:PT}  we may conclude:
\begin{prop}  Let $\exp_q^{nh}$ be the nonholonomic exponential map of a purely kinetic $(L,\D)$ and let $\hat{g}$ be a Riemannian metric which is a $\D$-preserving modification of $g$.  If, for all $v_q\in\D_q$, they  satisfy both 
$$
\hat{g}_{\exp^{nh}_q(v_q)}\left(T_{v_q}\left(\exp _q^{nh}\right)\left(v_{q}\right), T_{v_q}\left(\exp _q^{nh}\right)\left(w_{q}\right)\right)=\hat{g}_q(v_{q},w_{q}), \qquad \forall u_q \in \D_q 
$$
and   
$$ 
\hat{g}_{\exp_q^{nh}(tv_q)}\left(T_{tv_q}(\exp^{nh}_q)(tv_q),W(t)\right)+t \Big(h(t)\Big)\Big(W(t)\Big)=\hat{g}_q(tv_q,w), \qquad \forall w_q\in\D^g_q,
$$
 then ${\hat g}$ is a geodesic extension of $(L,\D)$.
\end{prop}
Because of the similarity with the statement of the Gauss Lemma, we call the conditions $(A)$ and $(B)$ the infinitesimal version of the above Gauss-type conditions.

\section{Chaplygin systems} \label{sec7}

Chaplygin systems (also called generalized Chaplygin systems or nonholonomic systems of {principal} type, see e.g.\ \cite{Bloch,Frans,LuisJC}) are essentially nonholonomic systems with extra symmetry properties. The goal of this section is to use that symmetry to simplify the conditions $(A)$ and $(B)$. 

Suppose that we have a purely kinetic Chaplygin system. This means that both the Riemannian metric $g$ (and its Lagrangian $L$) and the constraint manifold $\C$ are invariant under the action $G\times Q \to Q$ of a Lie group $G$ on the configuration manifold $Q$ (when lifted to $TQ$). In addition, we  assume that $\pi: Q\to Q/G$ is a principal fibre bundle and that the nonholonomic distribution $\mathcal{D}$ is  the horizontal distribution of a principal connection $\omega$ on $\pi$. Let $\{E_i\}$ be a basis for the Lie algebra of $G$ and $\{(E_i)_Q\}$ its corresponding infinitesimal generators. If $(x^a)$ are coordinates on $Q/G$, the horizontal lifts $X_a$ of the vector fields $\partial/\partial x^a$ on $Q/G$ (by means of $\omega$) are invariant vector fields on $Q$ that span $\D$. 
 The set $\{X_a,(E_i)_Q\}$ is then a basis for $\vectorfields{Q}$.
If we construct the basis in this way, we get
$$
[X_a, (E_i)_Q]=0,\quad[X_a, X_b]=B_{a b}^j (E_j)_Q,\quad[(E_i)_Q, (E_j)_Q]=-C_{i j}^k (E_k)_Q.
$$
Here, $B_{a b}^j$ can be interpreted as the curvature of the connection $\omega$ and $C_{i j}^k$ as the structure constants of the Lie algebra.

Let, for now, $G_{ab} = g(X_a,X_b)$, $G_{ai} = g(X_a, ({E}_i)_Q)$ and $G_{ij} = g(({E}_i)_Q,({E}_j)_Q)$. Similarly as  in Section~\ref{sec2}, we can create a basis  for $\mathcal{D}^g$, by defining the vector fields ${\tilde X}_i$ as 
$$
{\tilde X}_i=({E}_i)_Q- {G}^{ab}{G}_{ai}X_b=({E}_i)_Q+ {K}^{b}_{i}X_b.
$$ 
In the orthogonal basis $\{X_a,{\tilde X}_i\}$, we have  $g_{ab} = g(X_a,X_b)=G_{ab}$, $g_{aj} = g(X_a,\tilde{X}_j)=0$ and $g_{ij} = g(\tilde{X}_i,\tilde{X}_j)$.  We first compute the  bracket coefficients $R^k_{ia}$ and $R^k_{bc}$  in terms of the current bracket coefficients $B^k_{ab}$ and $C^k_{ij}$. 

On the one hand, we have
$
[X_a, {\tilde X}_i]= [X_a, ({E}_i)_Q+ {K}^{b}_{i} X_b ]
= X_a ({K}^{b}_{i}) X_b+{K}^{b}_{i} B_{a b}^j ({E}_j)_Q$, while on the other hand is $[X_a, {\tilde X}_i ]=R_{a i}^b X_b+R_{a i}^j {\tilde X}_j=  (R_{a i}^b+ {K}^{b}_{j} R_{a i}^j) X_b+R_{a i}^j ({E}_j)_Q$. Therefore,  $R_{a i}^j= {K}^{b}_{i} B_{a b}^j$. 
Likewise, from $[X_a, X_b]=B_{a b}^j ({E}_j)_Q$ and $[X_a, X_b]=R_{a b}^c X_c+R_{a b}^j {\tilde X}_j  = (R_{a b}^c-K^c_j R_{a b}^j) X_c+R_{a b}^j ({E}_j)_Q$, we can conclude that $R_{a b}^j= B_{a b}^j$.

The Lagrangian is still $L= \frac{1}{2}g_{ab}v^a v^b+\frac{1}{2}g_{ij}v^i v^j$. We now  give an expression for the multipliers $\lambda_i$ in the Chaplygin case. From the invariance of the Lagrangian it follows that $({E}_i)_Q^C(L)=0$ (see e.g.\ \cite{Inv}). Moreover, on $\C$, $X_b^C(L)=\Gamma_{(L, D)}\left(X_b^V(L)\right)$.  With that, we  have, on $\C$,
$$
\begin{aligned}
\lambda_i & =\Gamma_{(L, \D)}\left({\tilde X}_i^V L\right)-{\tilde X}_i^C(L)   =\Gamma_{(L, \D)}\left(g_{i j} v^j\right)- {K}^{b}_{i} X_b^C(L)- {\dot{K}^{b}_{i}} X_b^V(L) \\
& =-{K}^{b}_{i} \GLD \left(X_b^V(L)\right)- {\dot{K}^{b}_{i}} X_b^V(L)  =-\Gamma_{(L, \D)}\left({K}^{b}_{i} X_b^V(L)\right) =
\Gamma_{(L,\mathcal{D})}\left(G_{ai} v^a\right).
\end{aligned}
$$

We may now re-express the conditions $(A)$ and $(B)$. 
The $(A)$ condition, $\theta_k R^k_{ac} v^a =0$, becomes here $\theta_k B^k_{ac} v^a =0$. For the $(B)$ condition, we obtain after plugging in $(A)$ and the expression for $\lambda_i$,
\begin{eqnarray*}
0&=&\Gamma_{(L,\mathcal{D})}(\theta_i)+\lambda_i+\theta_k R^k_{ia}v^a =   \Gamma_{(L,\mathcal{D})}(\hat{g}_{ai}v^a)+ \Gamma_{(L,\mathcal{D})}\left(G_{ai} v^a\right) - \theta_k    K^{b}_{i} B_{a b}^k  v^a\\
  &=& \Gamma_{(L,\mathcal{D})}\left((\hat{g}_{ai}+G_{ai}) v^a\right).
\end{eqnarray*}	

As a consequence, $\mu_i = (\hat{g}_{ai}+G_{ai}) v^a$ are first integrals of the nonholonomic vector field $\GLD$.

Recall that $(G_{ai})$ in the Proposition below stands for the restriction $g|_{\D \times V\pi}$, where $V\pi$ is the vertical distribution of the principal fibre bundle $\pi: Q \to Q/G$.

\begin{prop} \label{propchap} Let $(L,\D)$ be a purely kinetic Chaplygin nonholonomic system. If $\GLD$ has  linear first integrals  $\mu_i=\mu_{ai}v^a$ that satisfy
\[
\mu_k B^k_{ac} v^a = G_{bk} B^k_{ac} v^bv^a,
\] 
then there exists a geodesic extension of $(L,\D)$ by a Riemannian metric $\hat g$.
 \end{prop}

\begin{proof}
{From the discussion above we may deduce that the condition $(B)$ is equivalent with asking for the existence of  first integrals  $\mu_i$ of $\GLD$. If these satisfy the conditions of the current proposition, $\theta_i=\mu_i - G_{ai} v^a$ will satisfy condition $(A)$.
Therefore,} the proof follows from Proposition~\ref{prop:result2}, since conditions $(A)$ and $(B)$ ensure that a (0,2) tensor field ${\hat g}$ with ${\hat g}(X_a,X_b) =g_{ab}$ and ${\hat g}(X_a,{\tilde X}_i) =\mu_{ai}-G_{ai}$ can be extended to a Riemannian metric.
\end{proof}

 Also in this approach we can derive further necessary algebraic conditions, similar to the one we had stated in Proposition~\ref{propsolv} for $\theta_k$. Recall from the above that the condition $(A)$ can be written as  $\theta_k B^k_{ac}v^a=0$ and that the condition $(B)$ is here $\GLD(\theta_k+G_{bk}v^b)=0$. When we apply the vector field $\GLD$ on the $(A)$-condition, we get after substituting $(B)$,
 \[
(C) \qquad \quad    \theta_k \GLD({S}^k_c)= - \lambda_k {S}^k_c,
		\]
where ${\lambda}_k=\GLD(G_{bk}v^b)$ and ${S}^k_c=R_{ac}^k v^a = B_{ac}^k v^a$. This new algebraic $(C)$-condition can be used, together with the algebraic condition $(A)$, to narrow down the possible solutions $\theta_k$ of the geodesic extension problem. We will apply this technique on the example of the two-wheeled carriage in Section~\ref{carriagesec}.

\section{Examples}
 
\subsection{The vertically rolling disk}

 The purpose of this section is to determine a geodesic extension  $\hat{g}$, on a concrete example. We have chosen the example of the disk (rolling vertically without slipping), not only  because it is well-known and simple, but also  because it is involved enough to enlighten some of the constructions and results we have obtained in the previous sections. Moreover, for this example, already, some quadratic Lagrangians $\hat L$ are known in the literature \cite{BFM,Simoes} and at the end of the section we will compare them with the ones we find through our methods.

We will use the following coordinates, see the figures in \cite{FB} or in \cite{TomBrussel}. The  pair $(x,y)$ stands for the position of the center of the mass $C$, 
$\varphi$ is the orientation angle of the disk (measured between the tangent of the contact point and the positive $x$-axis) and  $\theta$ is the  angle  that a point fixed
on the disk makes with respect to the vertical. In this example the configuration manifold is $Q=\mathbb{R}^2\times S^1\times S^1$. Because we assume the disk to remain vertical during its motion, the potential is constant and the Lagrangian is of purely kinetic type. The Lagrangian function is given by 
$$
L=\frac{1}{2}g(v,v) = \frac{1}{2}(\m\dot{x}^2+\m\dot{y}^2+\I\dot{\theta}^2+\J\dot{\varphi}^2), \qquad \I = \frac{1}{2}\m\Ra^2, \J = \frac{1}{4}\m\Ra^2. 
$$
and the linear nonholonomic constraints are 
\begin{equation}
\begin{cases}
\dot{x}=\Ra\cos(\varphi)\dot{\theta},\\
\dot{y}=\Ra\sin(\varphi)\dot{\theta}.
\end{cases} \nonumber
\end{equation}
This means that the distribution $\mathcal{D}$ can be spanned by 
$$
\mathcal{D}=\text{span}\{X_a\}=\textrm{span}\left\{X_\theta=\frac{\partial}{\partial \theta}+\Ra\cos(\varphi)\frac{\partial}{\partial x}+\Ra \sin(\varphi)\frac{\partial}{\partial y}, X_\varphi=\frac{\partial}{\partial \varphi}  \right\}.
$$ 
The constrained Lagrangian is 
\[
L_c = L|_\C =  \frac{1}{2}(\J\dot{\varphi}^2+(\m\Ra^2+\I)\dot{\theta}^2). 
\]

With this input one may write down the Lagrange-d'Alembert equations. It is well-known (see e.g.\ \cite{Bloch,Cortes,BFM}) that the center of mass will make either a circular movement, or a motion along a straight line.

\subsubsection{The vertically rolling disk as a purely kinetic nonholonomic system} \label{sec81}

Remark first that the metric $g$ is such that 
$$
g({\partial}/{\partial q^i},{\partial}{/\partial q^j})= \begin{pmatrix} \m & 0& 0& 0\\0 & \m &0& 0 \\ 0&0&\I&0 \\0&0&0&\J  \end{pmatrix}
$$
with $(q^i)=(x,y,\theta,\varphi)$. 

 In this subsection we choose the following basis for the orthogonal complement $\mathcal{D}^g$:
$$
\mathcal{D}^g=\text{span}\{X_i\}=\text{span}\left\{X_x=\frac{\partial}{\partial x}-\ce\cos(\varphi)\frac{\partial}{\partial \theta},X_y=\frac{\partial}{\partial y}-\ce\sin(\varphi)\frac{\partial}{\partial \theta}\right\},
$$
with $\ce = \frac{\Ra\m}{\I} = \frac{2}{R}$. It is important to realize that none of our results depend on the specific choices we have made for the frames. We will make this clear in the next subsection.

When compared to the expressions in the previous sections, the index $a$ will run over $\varphi$ and $\theta$, while the index $i$ is here either $x$ or $y$. If we compute the metric coefficients with respect to this basis, we get $g_{ai}=0$, meaning here that $g_{\theta x} = g_{\theta y} = g_{\varphi x}= g_{\varphi y}= 0$. Moreover 
\[
(g_{ab}) = \begin{pmatrix} g_{\theta \theta} & g_{\theta\varphi} \\ g_{\theta\varphi} & g_{\varphi \varphi} \end{pmatrix}  = \begin{pmatrix} \I+\m\Ra^2     & 0 \\ 0 & \J    \end{pmatrix} = MR^2\begin{pmatrix}  \frac{3}{2}    & 0 \\ 0 &     \frac{1}{4}   \end{pmatrix}
\]
and
\begin{eqnarray*}
(g_{ij}) &=& \begin{pmatrix} g_{xx} & g_{xy} \\ g_{xy} & g_{yy} \end{pmatrix} = \begin{pmatrix} \m+\ce^2\I\cos^2(\varphi) & \ce^2\I\cos(\varphi)\sin(\varphi) \\ \ce^2\I\cos(\varphi)\sin(\varphi)  & \m+\ce^2\I\sin^2(\varphi) \end{pmatrix} \\&=&  M\begin{pmatrix} 1+2 \cos^2(\varphi) & 2 \cos(\varphi)\sin(\varphi) \\ 2 \cos(\varphi)\sin(\varphi)  & 1+2 \sin^2(\varphi) \end{pmatrix}.
\end{eqnarray*}
With this, the Lagrangian in quasi-velocities $(v_x, v_y, v_\theta, v_\varphi)$ becomes
\begin{eqnarray*}
L&=&\frac{1}{2}\left(\left(\m+\ce^2\I\cos^2(\varphi)\right) v_x^2+2 \ce^2\I\cos (\varphi) \sin(\varphi) v_x v_y\right. 
\left.+\left(\m+\ce^2\I\sin^2(\varphi)\right) v_y^2\right) \\&& +\frac{1}{2}\left((\I+\m\Ra^2) v_\theta^2+\J v_{\varphi}^2\right) 
\\&=&
\frac{M}{2}\left(\left(1+2 \cos^2(\varphi)\right) v_x^2+4 \cos (\varphi) \sin(\varphi) v_x v_y\right. 
\left.+\left(1+2 \sin^2(\varphi)\right) v_y^2\right)\\&&+ \frac{MR^2}{2}\left(\frac{3}{2} v_\theta^2+\frac{1}{4} v_{\varphi}^2\right).
\end{eqnarray*}

We may read off the  coefficients $R_{\alpha \beta}^\gamma$ from the Lie brackets $\left[X_\alpha, X_\beta\right]=R_{\alpha \beta}^c X_c+R_{\alpha \beta}^k X_k$. Here they are
\begin{eqnarray*}
&& \left[X_{\varphi}, X_\theta\right]=-\Ra\sin(\varphi) X_x+\Ra\cos(\varphi) X_y,\qquad  \left[X_x, X_y\right]=0, \qquad
\left[X_x, X_\theta\right]=0, \qquad 
\left[X_y, X_\theta\right]=0,
\\ &&
\left[X_{\varphi}, X_x\right]=
\ct \sin(\varphi) X_\theta-\cd \sin(\varphi)\cos(\varphi) X_x-\cd \sin^2(\varphi) X_y 
,\\ && 
\left[X_\varphi, X_y\right] = -\ct \cos(\varphi) X_\theta+\cd \cos^2(\varphi) X_x+\cd \cos(\varphi)\sin(\varphi) X_y,
\end{eqnarray*}
with $\ct=\frac{\m\Ra}{\I +\m\Ra^2}= \frac{2}{3R}$ and $\cd =\frac{\m\Ra^2}{\I +\m\Ra^2} =  \frac{2}{3}$.

Since all $g_{ab}$ are constant and $R^a_{db}=0$, we get from the 
Koszul formula,
\[
2 g_{c d} \Gamma_{a b}^d =  X_a\left(g_{b c}\right)+X_b\left(g_{a c}\right)-X_c\left(g_{a c}\right)+g_{d c} R_{b c}^d
 -g_{d b} R_{a c}^d-g_{a d} R_{b c}^d,
\]
 that all $\Gamma_{a b}^d=0$.

We are now ready to find a solution for the equations $(A)$ and $(B)$ of Lemma~\ref{prop:result}.
For $(A)$, we get
the following two algebraic conditions that should be satisfied by $\hat{g}_{\varphi x},\hat{g}_{\varphi y},\hat{g}_{\theta x}$ and $\hat{g}_{\theta y}$: 
$$ 
\begin{cases}
\hat{g}_{\varphi x} \sin (\varphi)-\hat{g}_{\varphi y} \cos (\varphi)=0, \\
 -\hat{g}_{\theta x} \sin( \varphi)+\hat{g}_{\theta y} \cos (\varphi)=0.  
\end{cases}
$$

Condition $(B)$, written in a version without $v^av^b$  (and with $\Gamma^d_{bc}=0$) is here: {$$
 0=X_a(\hat{g}_{bi}) + X_b(\hat{g}_{ai})-X_i(g_{ab})-g_{db}R^d_{ai} -\hat{g}_{bk} R^k_{ai}-g_{ad}R^d_{bi}  -\hat{g}_{ak} R^k_{bi}.
$$}
If we first take $a=b=\varphi$, this gives $X_{\varphi}\left(\hat{g}_{\varphi_i}\right)-\hat{g}_{\varphi k} R_{\varphi_i}^k=0$, or:
\[
X_{\varphi}\left(\hat{g}_{\varphi x}\right)-\hat{g}_{\varphi x} R_{\varphi x}^x-\hat{g}_{\varphi y} R_{\varphi x}^y=0 \quad\mbox{and}\quad X_{\varphi}\left(\hat{g}_{\varphi y}\right)-\hat{g}_{\varphi x} R_{\varphi y}^x-\hat{g}_{\varphi y} R_{\varphi y}^y=0.
\]
After substituting the values of $R_{\varphi x}^y$, etc, and after substituting the relations we had derived from the algebraic condition  $(A)$, we may write this as
\[
\frac{\partial \hat{g}_{\varphi x}}{\partial \varphi}+\cd \hat{g}_{\varphi x} \tan (\varphi)=0
\quad \mbox{and}\quad \frac{\partial \hat{g}_{\varphi y}}{\partial \varphi}-\cd \hat{g}_{\varphi y} \cot( \varphi)=0. 
\]
The solutions of these equations are $\hat{g}_{\varphi x}=C|\cos (\varphi)|^{\cd}$  and $\hat{g}_{\varphi y}=D |\sin(\varphi)|^{\cd}$ (with $C,D$ possibly functions of $(x,y,\theta)$), but the only ones satisfying again the algebraic conditions are $\hat{g}_{\varphi x}=\hat{g}_{\varphi y}=0$.

When we put $a=\theta$ and $b=\varphi$, we get:
\[
0=X_{\varphi}\left(\hat{g}_{\theta i}\right)-g_{\theta\theta} R_{\varphi i}^\theta-\hat{g}_{\theta k} R_{\varphi i}^k,
\]
which becomes here:
\[
\frac{\partial \hat{g}_{\theta x}}{\partial \varphi}-g_{\theta\theta}\ct\sin (\varphi)+\cd \hat{g}_{\theta x} \sin( \varphi) \cos( \varphi)+\cd {\hat g}_{\theta y} \sin^2( \varphi)=0
\]
and
\[
\frac{\partial \hat{g}_{\theta y}}{\partial \varphi}+g_{\theta\theta}\ct\cos( \varphi)-\cd \hat{g}_{\theta x} \cos^2 (\varphi)-\cd {\hat g}_{\theta y} \cos(\varphi) \sin(\varphi)=0,
\]
from which (with $E,F$ functions of $(x,y,\theta)$)
\[
\hat{g}_{\theta x} =E|\cos (\varphi)|^{\cd}+\cv\cos( \varphi) \quad\mbox{and}\quad 
\hat{g}_{\theta y} = F|\sin (\varphi)|^{\cd}+\cv\sin( \varphi), 
\]
with $\cv = -\frac{\m\Ra}{\I}(\I+\m\Ra^2) = -3MR$.  Again, only for $E=0$ and $F=0$ the algebraic conditions are satisfied.  Finally, we get 
for $a=b=\theta$ that $X_\theta\left(\hat{g}_{\theta i}\right)=0$ (since $R_{\theta i}^\alpha=0$). This is automatically satisfied since the $\hat{g}_{\theta i}$ we had found only depend on $\varphi$.

We conclude that the trajectories of the vertically rolling disk with initial velocities in $\mathcal{D}$ are geodesics of all metrics $\hat{g}$ with metric coefficients in the frame $\{X_a,X_i\}$ given by 
$$
\hat{g}_{ab}=g_{ab},\quad\quad \hat{g}_{\theta x}=\cv \cos (\varphi), \quad \quad \hat{g}_{\varphi x}=0, \quad\quad \hat{g}_{\theta y}=\cv \sin( \varphi), \quad \quad\hat{g}_{\varphi y}=0,
$$
 and where the coefficients $\hat{g}_{ik}$ can be anything as long as $\hat{g}$ is positive-definite. 
We may look at the constructive proof of Proposition \ref{prop:result2} for one possible (but not exclusive) way  to realize this: we only need to find a function $\alpha$ such that $\hat{g}=\begin{pmatrix}
    A&B\\B^T&\alpha I
\end{pmatrix}$ is positive-definite. The matrices $A$ and $B$ are here given by 
$$A=\begin{pmatrix}
    \I+\m\Ra^2&0\\0&\J
\end{pmatrix}\quad\text{ and }\quad B=\begin{pmatrix}
    \cv\cos(\varphi)&\cv\sin(\varphi)\\0&0
\end{pmatrix}.
$$
 When we take  $\alpha$ to be a constant, it only needs to be greater than the  largest eigenvalue of $B^TA^{-1}B$. Here, these eigenvalues are $0$ and $\frac{\cv^2}{I+\m\Ra^2}= 6M$. We see that any $\alpha>6M$ fulfills all requirements.

If we set, for simplicity, $M=1$ and $R=1$, we  can take for example $\alpha=42$ and then 
$$
\hat L =\frac{1}{2}\left[\frac{3}{2}v_\theta^2 +\frac{1}{4}v_\varphi^2-6\cos(\varphi) v_\theta v_x-6\sin(\varphi) v_\theta v_y+42v_x^2+42v_y^2\right].
$$

The relation between the quasi-velocities $\{v_x, v_y, v_\theta, v_\varphi\}$ and the natural fibre bundle coordinates   $\{ \dot{x},\dot{y},\dot{\varphi},\dot{\theta}\}$ is
\[
\dot{x}=v_x+v_\theta \cos( \varphi),\quad \dot{y}=v_y+v_\theta \sin(\varphi),\quad  \dot{\varphi}=v_\varphi, \quad \dot{\theta}=v_\theta-  2 v_x \cos (\varphi)-  2 v_y \sin (\varphi),
\]
or conversely  
\[ \begin{cases}
v_\varphi=\dot\varphi,\\ v_\theta = \frac{1}{ 3 }(\dot\theta +  2 \cos(\varphi)\dot x +  2 \sin(\varphi)\dot y),\\ v_x = \frac{1}{ 3 }(-\cos(\varphi)\dot\theta + ( 2 \sin^2(\varphi)+1)\dot x -  2 \cos(\varphi)\sin(\varphi)\dot y), \\  v_y = \frac{1}{ 3 }(-\sin(\varphi)\dot\theta -  2 \cos(\varphi)\sin(\varphi)\dot x +  (  2 \cos^2(\varphi)+1)\dot y).
\end{cases}
\]
After substituting this, one obtains an expression for $\hat L$ in natural coordinates.

\begin{center}
\begin{tabular}{cr}
\begin{minipage}{6.5cm}

\includegraphics[scale=0.3]{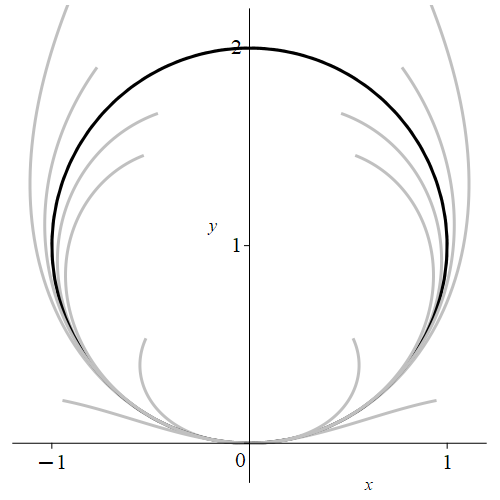}
\end{minipage}
&
\begin{minipage}{8.5cm}

On the left, we have plotted some  geodesics  for $\hat g$ with $\alpha=42$ (their projections for the center of mass $(x,y)$). The black circle through the origin  is the nonholonomic trajectory through the origin with initial velocity \[
({\dot x}_0=1, {\dot y}_0=0,  {\dot \theta}_0=1, {\dot \varphi}_0=1)\in \D_{(0,0,0,0)}.
\] The nearby grey geodesics are perturbations in either ${\dot x}_0$, ${\dot y}_0$ and ${\dot \theta}_0$, with $\epsilon = \pm 0.01$.
\end{minipage}
\end{tabular}
\end{center}

We can now compare our results with the Lagrangians that have been derived in the literature, as a solution of the Hamiltonization problem. A second attempt to obtain a Riemannian metric $\hat g$ from the same solution $(g_{ab},{\hat g}_{ai})$ is to choose ${\hat g}_{ij}$ in such a way that it is a multiple of the original $g_{ij}$, i.e.\ ${\hat g}_{ij} = \beta g_{ij}$. 
In \cite{BFM} (see also \cite{FB}) we can find the following regular quadratic Lagrangian (among other non-quadratic Lagrangians,   and after setting $M=1$ and $R=1$):
\[
{\hat L}_{\cite{BFM}} = \frac{1}{2} ( \frac{1}{2} {\dot\theta}^2 +  \frac{1}{4}  {\dot\varphi}^2) - \frac{1}{2}( {\dot x}^2 + {\dot y}^2) + \dot\theta (\dot x \cos(\varphi) + \dot y \sin(\varphi)).
\]
When written in our current quasivelocities, we find for the corresponding metric 
\[
({\hat g}_{ab})_{\cite{BFM}} = ({\hat g}_{ab}),\quad ({\hat g}_{ai})_{\cite{BFM}} = ({\hat g}_{ai}),\quad ({\hat g}_{ij})_{\cite{BFM}} = - \begin{pmatrix} 1+ {2} \cos^2(\varphi) &   {2} \cos(\varphi)\sin(\varphi) \\[2mm]  {2} \cos(\varphi)\sin(\varphi) & 1+ {2} \sin^2(\varphi)  \end{pmatrix},
\]
meaning that it corresponds with $\beta=-1$. This metric has therefore the same constrained Lagrangian $L_c$ as $L$ and as our Lagrangian $\hat L$. However, the choice $\beta=-1$ leads only to a pseudo-Riemannian metric.

The following Lagrangian is mentioned in \cite{Simoes}:
\[
{\hat L}_{\cite{Simoes}} = \frac{1}{2} (2{\dot\theta}^2 + {\dot\varphi}^2+ {\dot x}^2 + {\dot y}^2 ) - \dot\theta (\dot x \cos(\varphi) + \dot y \sin(\varphi)). 
\]
It  comes from a Riemannian metric, but it does not preserve the constrained Lagrangian (i.e.\ it is not a $\D$-preserving modification of $g$), since 
\[
(g_{ab})_{\cite{Simoes}} =  \begin{pmatrix} 1 & 0 \\ 0 & 1 \end{pmatrix}.
\]

 The metric that is denoted by $h$ in \cite{recent} (in the example section on the vertically rolling disk) is a $\D$-preserving modification, and as such it belongs to our class of metrics.

 Finally, there are also some results that can be related to the geodesic extension problem in \cite{BC}. However, the authors of \cite{BC} do not make use of a $\D$-preserving modification, and they mention that they are unable to find a geodesic extension by their methods, specifically, for the vertically rolling disk.

\subsubsection{The vertically rolling disk as a Chaplygin system}

It is well-known (\cite{Bloch,Cortes}) that the configuration space $Q=S^1\times S^1\times \mathbb{R}^2$ of the vertically rolling disk is the total space of a principal fibre bundle with structure group $G=\mathbb{R}^2$. Both the Lagrangian and the constraints are invariant under this symmetry group, and the constraints can be thought of as the horizontal distribution of a principal connection. This means that the vertically rolling disk can be regarded as a Chaplygin system.

We consider again the case with $M=R=1$.  For the vertical distribution of $\pi:Q\rightarrow Q/G$ we can take $V\pi=\text{span} \{({E}_x)_Q=\frac{\partial}{\partial x} ,({E}_y)_Q=\frac{\partial }{\partial y}\}$, and for the horizontal distribution we may set $\mathcal{D}=\text{span}\{X_\varphi=\frac{\partial}{\partial \varphi},X_\theta=\cos(\varphi)\frac{\partial}{\partial x}+ \sin(\varphi)\frac{\partial}{\partial y} +\frac{\partial}{\partial \theta} \}$, as before. We may compute that  
\[
G_{x \varphi}=0, \quad
G_{x \theta}=\cos( \varphi), \quad G_{y \varphi}=0, \quad G_{y \theta}=\sin( \varphi),
\]
and
\[ \left[X_{\varphi}, X_\theta\right]=-\sin (\varphi) \frac{\partial}{\partial x}+\cos (\varphi) \frac{\partial}{\partial y}.
\]

It is important to remark, however, that the result in Proposition~\ref{propchap} makes use of a frame of $\D^g$ that is different from the frame $\{X_x,X_y\}$ we had encountered earlier. The vector fields 
${\tilde X}_i =({E}_i)_Q-G^{bc}G_{ci} X_b=({E}_i)_Q+K^b_i X_b$ are here 
$$
{\tilde X}_x=\frac{\partial}{\partial x}- \frac{2}{3} \cos(\varphi)\left(\frac{\partial}{\partial \theta}+\cos(\varphi)\frac{\partial}{\partial x}+\sin(\varphi)\frac{\partial}{\partial y}\right)
$$
 and
$$
{\tilde X}_y=\frac{\partial}{\partial y}- \frac{2}{3} \sin(\varphi)\left(\frac{\partial}{\partial \theta}+\cos(\varphi)\frac{\partial}{\partial x}+\sin(\varphi)\frac{\partial}{\partial y}\right),
$$
and the  relation between $\{{\tilde X}_i\}$ and the frame $\{{X}_i\}$ of the previous subsection is 
\begin{equation*}
\begin{cases}
    {\tilde X}_x=(1- \frac{2}{3} \cos^2(\varphi))X_x- \frac{2}{3} \cos(\varphi)\sin(\varphi) X_y,\\
    {\tilde X}_y=- \frac{2}{3} \sin(\varphi)\cos(\varphi) X_x+(1- \frac{2}{3} \sin^2(\varphi))X_y.
\end{cases}\nonumber
\end{equation*}

According to Proposition~\ref{propchap} we should look for linear first integrals $\mu_x$ and $\mu_y$  of $\Gamma_{(L,\mathcal{D})}$  that satisfy 
$$
 \left(\mu_x-G_{x \theta} v_\theta-G_{x \varphi} v_{\varphi}\right) B_{\varphi \theta}^x v_{\varphi} +\left(\mu_y-G_{y \theta} v_\theta-G_{y \varphi} v_\phi\right) B_{\varphi \theta}^y v_{\varphi}=0,
$$
which is here
$$
-\sin( \varphi) \mu_x+\cos (\varphi) \mu_y=0. 
$$
If we take a $\GLD$-derivative we also get $\cos( \varphi) \mu_x+\sin (\varphi) \mu_y=0$ (because $\GLD(\sin(\varphi)) = \cos(\varphi)\dot\varphi$). 
We can, therefore conclude that the first integrals should be $\mu_x=\mu_y=0$. With that, we may get the coefficients of the metric $\hat g$ with respect to the frame $\{X_a,{\tilde X}_i\}$ (see the proof of Proposition~\ref{propchap}),
\[
{\hat g}(X_a,X_b) =g_{ab}, \qquad \widetilde{{\hat g}_{ai}}= {\hat g}(X_a,{\tilde X}_i) =\mu_{ai}-G_{ai},
\]
which are here
\[    \begin{cases}
    \widetilde{\hat{g}_{\varphi x}}=0 \qquad \mbox{and}\qquad \widetilde{\hat{g}_{\theta x}}=-\cos( \varphi)\\
 \widetilde{\hat{g}_{\varphi y}}=0 \qquad \mbox{and}\qquad \widetilde{\hat{g}_{\theta y}}=-\sin( \varphi)
    \end{cases}.
		\]
We see that this method is way shorter and more efficient since we had already solved the involved differential equations. 

We end this section by showing that the new coefficients we had found are compatible with the old ones. Indeed, given the relation between the two bases $\{X_i\}$ and $\{{\tilde X}_i\}$ one easily verifies that, for example,  
\[
    \widetilde{\hat{g}_{\theta y}} =- \frac{2}{3}  \sin(\varphi) \cos(\varphi) \hat{g}_{\theta x}+\left(1- \frac{2}{3} \sin ^2 (\varphi)\right) \hat{g}_{\theta y}, 
\]
as it should.

\subsection{A two-wheeled carriage}  \label{carriagesec}

In this section, we will analyze the solvability of the conditions $(A)$ and $(B)$ for the example of a two-wheeled carriage (see also \cite{anhol,deL}). Since it is a Chaplygin system, we will first fully exploit the algebraic conditions that we had derived after  Proposition~\ref{propchap}. In other words, we will  first find all $\theta_k$ that solve 
\begin{eqnarray*}
    (A) && \theta_k S^k_a=0,\\
    (C) && \theta_k \GLD(S^k_a)= \lambda_k S^k_a.   
    \end{eqnarray*}  
and, afterwards, verify if any of them also satisfy 
\begin{eqnarray*}
    (B) &&
 \Gamma_{(L,\mathcal{D})}(\theta_k+G_{bk}v^b)=0.
\end{eqnarray*}

\begin{center}
\begin{tabular}{lr}
\begin{minipage}{7cm}
\hspace{-0.7cm}\includegraphics[scale=0.8]{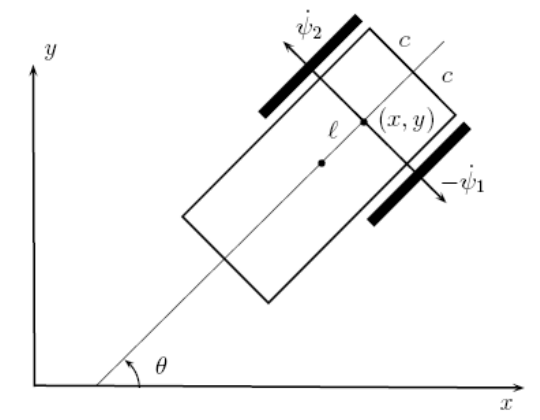}
\end{minipage}
& 
\begin{minipage}{8cm}
We will use  $\psi_1$ and $\psi_2$ for the angles that describe the rotation of the two wheels of the carriage, respectively. The coordinates $(x,y)$ describe the position of the intersection point of the horizontal symmetry axis of the carriage and the vertical axis that connects the two wheels. The angle $\theta$ is the angle between the symmetry axis and the $x$-axis.  The parameter $c$ stands for half of the distance between the two wheels. We will denote the   distance between the centre of mass  and  the  intersection point $(x,y)$  by the parameter $\ell$. We will think of  $\ell$ as a parameter that we can control when manufacturing the carriage.
\end{minipage}
\end{tabular}
\end{center}
The configuration space $Q$ is here $S^1 \times S^1 \times
SE(2)$. The Lagrangian function of the carriage is 
\[ L= \frac{1}{2} m ({\dot{x}}^2 + {\dot{y}}^2) +
m_0 \ell {\dot{\theta}}((\cos\theta) \dot{y} - (\sin\theta) \dot{x})
+\frac{1}{2} J {\dot{\theta}}^2 + \frac{1}{2} J_2 ({\dot\psi}_1^2 + {\dot\psi}_2^2 ),
\]
where $J_2$ and $J$ are the moments of inertia for the wheels and for the whole system, respectively, $m_0$ is the mass of the body and $m$ is the mass of the whole system. The two wheels of the carriage are assumed to roll without slipping and this can be described by the following three constraints:
\[
\dot{x} = -\frac{R}{2} \cos\theta (\dot{\psi_1} + \dot{\psi_2}),\quad
\dot{y} = -\frac{R}{2} \sin\theta (\dot{\psi_1} + \dot{\psi_2}), \quad
\dot{\theta} = \frac{R}{2c}  (\dot{\psi}_2 - \dot{\psi}_1).
\]
Here, $R$ is the radius of both wheels. From these constraints we can deduce a frame $\{X_{\psi_1}, X_{\psi_2}\}$ for the distribution $\mathcal{D}$: 
\begin{align*}
X_{\psi_1} &= \fpd{}{\psi_1} -
\frac{R}{2}\left(\cos\theta \fpd{}{x}+\sin\theta\fpd{}{y} +
\frac{1}{c}\fpd{}{\theta}\right), \\
X_{\psi_2} &= \fpd{}{\psi_2} -
\frac{R}{2}\left(\cos\theta \fpd{}{x}+\sin\theta\fpd{}{y} -
\frac{1}{c}\fpd{}{\theta}\right).
\end{align*}

The metric coefficients in the constrained Lagrangian $L_c=\frac{1}{2}g_{ab}v^a v^b$ are given by the matrix $$(g_{ab})=\begin{pmatrix}
   J_2+\frac{R^2m}{4}+\frac{R^2J}{4c}&\frac{R^2 m}{4}-\frac{R^2 J}{4c}\\[2mm]\frac{R^2 m}{4}-\frac{R^2 J}{4c}& J_2+\frac{R^2m}{4}+\frac{R^2J}{4c}
\end{pmatrix}=\begin{pmatrix}
    P&-Q\\[2mm]-Q&P
\end{pmatrix}.$$
Because of the Chaplygin properties, the  Lagrangian and the constraints are both invariant under the $SE(2)$-action. Therefore, we can consider the infinitesimal generators $(E_x)_Q,(E_y)_Q,(E_{\theta})_Q$ corresponding to a basis of the Lie algebra. They will span  the vertical distribution $V\pi$ of $\pi:Q\rightarrow Q/G$, where $Q=S^1 \times S^1 \times
SE(2)$ and $G=SE(2)$. These vector fields have the form 
\[
{E}_x=\fpd{}{x}, \quad
{E}_y=\fpd{}{y}, \quad
{E}_{\theta}=\fpd{}{\theta} - y \fpd{}{x} + x\fpd{}{y}.
\]

In this frame, the quasi-velocities $v^{\psi_1}$ and $v^{\psi_2}$ are (when restricted to $\D$) given by $\dot{\psi_{1}}$ and $\dot{\psi_2}$  and  we have 
\begin{align*}
\dot{x}&=-\frac{1}{2} R\cos\theta(v^{\psi_1}+v^{\psi_2}),\\
\dot{y}&=-\frac{1}{2} R\sin\theta(v^{\psi_1}+v^{\psi_2}),\\
\dot{\theta}&=-\frac{R}{2c}(v^{\psi_1}-v^{\psi_2}).
\end{align*}
The nonholonomic vector field $\Gamma_{(L,\mathcal{D})}$ in terms of the quasi-velocities $v^{\psi_1}$ and $v^{\psi_2}$ and the given frame   is of the form
$$\Gamma_{(L,\mathcal{D})}= v^{\psi_1}X_{\psi_1}^ C+v^{\psi_2}X_{\psi_2}^C+f^{\psi_1} X_{\psi_1}^V+f^{\psi_2} X_{\psi_2}^V,$$
where \[
f^{\psi_1}=\frac{K}{P^2-Q^2}(v^{\psi_1}-v^{\psi_2})(Qv^{\psi_1}-Pv^{\psi_2}),\quad
f^{\psi_2}=\frac{K}{P^2-Q^2}(v^{\psi_1}-v^{\psi_2})(Pv^{\psi_1}-Qv^{\psi_2}),
\]
with constant $K=m_0 \ell R^3/4c^2$.

Since 
\[
[X_{\psi_1},X_{\psi_2}]=\textcolor{red}{-}\frac{R^2}{2c}\left(\sin\theta(E_x)_Q-\cos\theta(E_y)_Q\right),
\]
the non-zero bracket coefficients $B^k_{ab}$  are given by \[
B^x_{\psi_1 \psi_2} =\textcolor{red}{-}\frac{R^2}{2c}\sin\theta \quad\mbox{and}\quad  B^y_{\psi_1 \psi_2} =\frac{R^2}{2c}\cos\theta.
\] 
Therefore, we get   for $S^k_a=B^k_{ab}v^b$,
$$S^x_{\psi_1}=\textcolor{red}{-}\frac{R^2}{2c}\sin\theta v^{\psi_2},\quad S^x_{\psi_2}=\frac{R^2}{2c}\sin\theta v^{\psi_1},\quad  S^y_{\psi_1}=\frac{R^2}{2c}\cos{\theta}v^{\psi_2},\quad S^y_{\psi_2}=\textcolor{red}{-}\frac{R^2}{2c}\cos\theta v^{\psi_1}.$$
From this and the expression of $\GLD$ one may easily compute the expressions for $\GLD(S^k_a)$. Likewise, the expressions for ${\lambda}_k=\GLD(G_{bk}v^b)$ follow from 
\begin{eqnarray*}
&& G_{\psi_1x} = -\frac{R}{2}m\cos\theta + \frac{R}{2c}m_0\ell \sin\theta,\qquad G_{\psi_1y} = -\frac{R}{2}m\sin\theta -\frac{R}{2c}m_0\ell \cos\theta,\\
&& G_{\psi_2x} = -\frac{R}{2}m\cos\theta - \frac{R}{2c}m_0\ell \sin\theta,\qquad G_{\psi_2y} = -\frac{R}{2}m\sin\theta +\frac{R}{2c}m_0\ell \cos\theta,
\end{eqnarray*}
 
For $a=\psi_1$ the algebraic condition $(A)$ gives us an expression for $\theta_y$ in terms of $\theta_x$:  
$$
  \displaystyle  \theta_x\left(\frac{R^2}{2c}\sin\theta v^{\psi_2}\right)-\theta_y\left(\frac{R^2}{2c}\cos\theta v^{\psi_2} \right)=0\quad\Leftrightarrow\quad \theta_y=\frac{\sin\theta}{\cos\theta}\theta_{\textcolor{red}{x}}.$$
 If we use this in the  algebraic condition $(C)$ for $a=\psi_1$, we get (after some computations) that the only potential values for $\theta_x$ and $\theta_y$ are 
  $$
	\theta_x=\left[ \frac{mR}{2}+\frac{K\ell m_0}{P+Q}\right]\cos\theta(v^{\psi_2}+v^{\psi_1}),\quad \theta_y= \left[ \frac{mR}{2}+\frac{K\ell m_0}{P+Q}\right]\sin\theta(v^{\psi_2}+v^{\psi_1}).
	$$
Since $\theta_x=\hat{g}_{x{\psi_1}}v^{\psi_1}+\hat{g}_{x{\psi_2}}v^{\psi_2}$ and $\theta_y=\hat{g}_{y{\psi_1}}v^{\psi_1}+\hat{g}_{y{\psi_2}}v^{\psi_2}$, we derive from this that $$\hat{g}_{x{\psi_1}}=\hat{g}_{x{\psi_2}}=\left[ \frac{mR}{2}+\frac{K\ell m_0}{P+Q}\right]\cos\theta\quad\textrm{and}\quad \hat{g}_{y{\psi_1}}=\hat{g}_{y{\psi_2}}=\left[ \frac{mR}{2}+\frac{K\ell m_0}{P+Q}\right]\sin\theta.  $$  
  The two algebraic conditions for $a=\psi_2$ lead to  the same expressions for $\theta_x$ and $\theta_y$.
	
Recall that the algebraic conditions $(A)$ and $(C)$ are only necessary, but not sufficient. We still need to check whether the two expressions above satisfy the PDE condition $(B)$, for both $k=x$ and $k=y$. That is, we need to check that 
$$\GLD(\theta_x+G_{x\psi_1}v^{\psi_1}+G_{x\psi_2 }v^{\psi_2})=0\quad\textrm{and}\quad \GLD(\theta_y+G_{y\psi_1}v^{\psi_1}+G_{y\psi_2 }v^{\psi_2})=0 .$$ 
When written in full, this gives  $$m_0 l\frac{(v^{\psi_1}-v^{\psi_2})^2}{4c^2(P^2-Q^2)}(4K^2c^2+(Q^2-P^2)R^2)\cos(\theta)=0 $$ for the first, and the same expression with $\sin(\theta)$ for the second. This means that, besides being zero, the distance $\ell$  to  the centre of mass  can only be such that $4K^2c^2+(Q^2-P^2)R^2=0$. After substituting the explicit expressions for  of $K,P$ and $Q$, we obtain two suitable values 
  $$ \ell=0\qquad \textrm{or} \qquad \ell=\frac{\sqrt{(R^2m+2J_2)(JR^2+2J_2c^2)}}{m_0 R^2}.$$
 Recall that the Lagrangian of a $\mathcal{D}$-preserving geodesic extension is given by $$\hat{L}=\frac{1}{2}(g_{ab}v^a v^b+2\hat{g}_{ai}v^a v^i+\hat{g}_{ij}v^iv^i),$$
where the metric coefficients $g_{ab}$ come from the original system, the $\hat{g}_{ai}$ are the ones we have computed and the coefficients $\hat{g}_{ij}$ may be chosen freely. When we plug in the  second value of $\ell$, this becomes $$\hat{L}=\frac{1}{2}\left[\left(\frac{R^2}{4c^2}(J+mc^2)+J_2\right)v^{\psi_1}v^{\psi_1}-2\frac{R^2}{4c^2}(J-mc^2) v^{\psi_1}v^{\psi_2} + \left(\frac{R^2}{4c^2}(J+mc^2)+J_2\right)v^{\psi_2}v^{\psi_2}\right]$$ $$+\left[(Rm+\frac{J_2}{R})\cos(\theta)v^xv^{\psi_1}+(Rm+\frac{J_2}{R})\cos(\theta)v^x v^{\psi_2} +(Rm+\frac{J_2}{R})\sin(\theta)v^y v^{\psi_1}+2(Rm+\frac{J_2}{R})\sin(\theta)v^y v^{\psi_2}\right]$$ $$+\frac{1}{2}\hat{g}_{ij}v^i v^j.$$ 

In conclusion, there are only two types of carriages that possess a $\D$-preserving geodesic extension. The second, rather surprising, value of $\ell$ is actually the same one that also appears in the context of a comparison of the nonholonomic and vakonomic dynamics in \cite{anhol}.

\section{Discussion and outlook}

In this paper we have derived, for kinetic nonholonomic systems, with or without a symmetry group, the conditions $(A)$ and $(B)$ for the existence of a geodesic extension by a  $\D$-preserving modification. Besides its intrinsic relation to the nonholonomic exponential map, we have shown how the approach can be successfully applied on concrete examples. Moreover,  the example of the two-wheeled carriage  shows that the conditions (together with further integrability conditions, such as $(C)$) can also be used, when one wishes to determine obstructions for the existence of a $\D$-preserving geodesic extension  (in this case on the value of $\ell$).  

One may find a second example of such an obstruction in the number of nonholonomic constraints. Indeed, we will  verify below that, generically, there will not exist a genuine ${\mathcal D}$-preserving modification when there is just a single  constraint.  In this case, the $(A)$-condition (in the unknowns  ${\hat g}_{1a}$) is ${\hat g}_{b1}R^1_{ac}+{\hat g}_{a1}R^1_{bc} = 0$. When we set $a=b$, we get ${\hat g}_{a1}R^1_{ac} = 0$ (no sum over $a$). Since $R^1_{ac} \neq 0$ (if not, the constraint is, in fact, holonomic), we obtain ${\hat g}_{a1}=0$. This means that we only have the freedom in ${\hat g}_{ij}$ to obtain a ${\mathcal D}$-preserving modification. In general that is not enough: Because of $(B)$, ${\hat g}_{a1}=0$ can only be a solution when the original system satisfies $g_{11}\Gamma^1_{ab} =0$, or when $\Gamma^1_{ab} =0$. 

We may test this on the nonholonomic particle with Lagrangian  \[
L = \frac{1}{2}({\dot x}^2+{\dot y}^2 +{\dot z}^2) + \alpha {\dot y}{\dot z}
\]
and  constraint  ${\dot z} = y{\dot x}$. In the orthogonal frame \[\left\{X_x = \fpd{}{x} + y \fpd{}{z}, X_y = \fpd{}{y}, X_1 = \fpd{}{z}+y^2(\alpha^2-1)\fpd{}{x} -\alpha \fpd{}{y} \right\},\]
one may verify that that
\[
\Gamma^1_{xx} = \frac{1}{D}(\alpha (3(1-\alpha^2)y^2 - 2 (1-\alpha^2) y +1) ), 
\]
(with $D$ some denominator of no concern). This is only zero when $\alpha=0$. But, if we now set $\alpha=0$,  one may also calculate that 
$\displaystyle 
\Gamma^1_{xy} = \frac{1-y}{2(1+y^2)(1+y)} \neq 0$,
which means that there does not exist an  $\alpha$ that gives a geodesic extension. 

In a future work we wish to overcome this obstruction by allowing a larger class of geodesic extensions than only those that are $\mathcal D$-preserving modifications.  For example, one could allow  for a conformal change ${\hat g}_{ab} =e^\phi g_{ab}$ and mitify the effect by a further parameter transformation along the geodesics. The case $\phi=0$ then corresponds with our $\D$-preserving modifications. In fact, this is the driving idea behind the recent preprint \cite{recent} and, actually also behind the so-called Chaplygin multiplier approach in the larger context of Hamiltonization, see e.g.\ \cite{Tom,luis}. In the paper \cite{recent} the authors were able to show that a geodesic extension exists after a conformal change, in the case where the nonholonomic system is Chaplygin with a  $\phi$-simple gyroscopic tensor (such as, e.g.,  the nonholonomic particle  with $\alpha=0$, above). But, in the overlap  $\phi=0$, our approach is more general in that we  do not assume the presence of a symmetry group and that we do not give just one  geodesic extension, but emphasize that there is a whole class of them, when one exploits the freedom in choosing the part in  ${\hat g}_{ij}$. Even if the system is Chaplygin, this ${\hat g}_{ij}$ can still be a part that does not reduce to the quotient space. 

The condition of $\phi$-simplicity on which \cite{recent} relies, is heavily related to the existence of an invariant measure (see e.g.\ \cite{LuisJC}).
The possible relation between geodesic extensions on the one hand and invariant measures and $\phi$-simplicity on the other hand, is an interesting path for future research.
Consider the bracket $[X_a,X_b] = R_{ab}^c X_c + R_{ab}^i X_i$. A Chaplygin nonholonomic system is said to be $\phi$-simple if the so-called gyroscopic tensor $R_{ab}^c$ has a prescribed form that can be derived from a single function $\phi$ on the reduced space.  However, if we were given a  nonholonomic system with an invariant measure, there is no obvious guess of how it could be related to the existence of a geodesic extension.  For example, for the special  $\ell\neq 0$ of  the two-wheeled carriage we did find a $\mathcal D$-preserving geodesic extension. Since this represents  a `conformal change with $\phi=0$'  its gyroscopic tensor should vanish for it to be   $\phi$-simple. This is not the case for our system: in fact, it is here
\[
R^{\psi_1}_{\psi_1\psi_2} = -R^{\psi_2}_{\psi_1\psi_2} = -\frac{R^3 m_0 }{4c^2(P+Q)} \ell.
\]

In a next paper we wish to generalize  the $(A)$ and $(B)$ conditions to allow for conformal changes. In contrast to \cite{recent}, we wish to consider also nonholonomic systems without symmetry groups, and we wish to drop the requirement of $\phi$-simplicity. Besides, we wish to clarify the role of an invariant measure in this larger context.

\paragraph{Acknowledgements.}  We are grateful to all referees. Without doubt, their remarks have improved our manuscript. We also  thank the   Research Fund of the University of Antwerp (BOF) for its support through the DOCPRO project  49747.

\end{document}